\documentclass[12pt]{article}
\renewcommand{\theequation}{\thesection.\arabic{equation}}

\newcommand{\IN}{I\!\!N}
 
\newcommand{\N}       {{\cal N}}
\renewcommand{\Re}{I\!\!R}

\newcommand{\ai}               {A_i}
\newcommand{\ais}               {A_i^*}
 
\newcommand{\bd}               {{\udelta}}
\newcommand{\bdi}               {\bar{\udelta}_i}
\newcommand{\bdis}               {\bar{\udelta}_i^*}
 
\newcommand{\cdi}               {\check{\udelta}_i}
\newcommand{\cdis}               {\check{\udelta}_i^*}

\newcommand{\ecd}               {E_{\cdot|\hat{\udelta}}}
\newcommand{\hub}               {\hat{\ub}}
\newcommand{\hud}               {\hat{\udelta}}
\newcommand{\hV}               {\hat{V}}
\newcommand{\moi}               {M_{1i}}
\newcommand{\mtoi}               {M_{2i}}
\newcommand{\uSi}               {{\uSigma}}

\newcommand{\beas}               {\begin{eqnarray*}}
\newcommand{\eeas}               {\end{eqnarray*}}
\newcommand{\De}               {{\Delta}}
\newcommand{\ep}               {{\epsilon}}
\newcommand{\mti}               {m\rightarrow\infty}
\newcommand{\raw}               {\rightarrow}

\newcommand{\avelno}               {N_{0}^{-1}\sum_{l=1}^{N_0}}
\newcommand{\z}       {{z_{m}}}

\def\qmq#1{\quad\mbox{#1}\quad}

\newcommand{\ua}       {\mbox{$a$}} 
 
\newcommand{\ub}       {\mbox{$b$}}

\newcommand{\ue}       {\mbox{$e$}}

\newcommand{\ut}       {\mbox{$t$}}

\newcommand{\uX}       {\mbox{$X$}}
\newcommand{\ux}       {\mbox{$x$}} 
\newcommand{\uY}       {\mbox{$Y$}}
\newcommand{\uy}       {\mbox{$y$}}

\newcommand{\udelta}            {\mbox{$\delta$}}

\newcommand{\uiota}             {\mbox{$\uiota$}}

\newcommand{\ulambda}           {\mbox{$\lambda$}}

\newcommand{\upsi}              {\mbox{$\psi$}}

\newcommand{\uDelta}            {\mbox{$\Delta$}}

\newcommand{\uSigma}            {\mbox{$\Sigma$}}

\newcommand{\be}               {\begin{equation}}
\newcommand{\ee}               {\end{equation}}
\newcommand{\bea}              {\begin{eqnarray}}
\newcommand{\eea}              {\end{eqnarray}}
\newcommand{\ba}               {\begin{array}}
\newcommand{\ea}               {\end{array}}
\newcommand{\nn}               {\nonumber}

\begin{document}

\begin{center}

\large
{\bf Resampling Based Empirical Prediction: An Application 
to Small Area Estimation}\\
%\normalsize
%\vspace*{2ex}
%
{\small \sc By Soumendra  N. Lahiri}$^1$, {\small\sc Tapabrata Maiti}$^1$, {\small\sc Myron Katzoff}$^2$,
{\small \sc And} {\small \sc Van Parsons}$^2$\\[.1in]

$^1${\small \it Department of Statistics, Iowa State University, Ames, IA 50011;
snlahiri,taps\{@iastate.edu\}
} 
~ $^2${\small \it   NCHS/CDC,
3311 Toledo Road,
Hyattsville, Maryland 20782; 
mjk5,vlp1\{@cdc.gov\}
}\\[.1in]
{\small{\sc Summary}}\\[.1in]
\end{center}
 Best linear unbiased prediction is well known for its 
wide range of applications including small area estimation. While the 
theory is well established for mixed linear models
 and under normality  of the  error and mixing distributions, 
the literature is sparse for nonlinear 
mixed models under  nonnormality of the error or of the  
mixing distributions. This 
article develops a resampling based unified approach
 for predicting mixed effects under a generalized mixed 
model set up. Second order accurate nonnegative estimators of 
mean squared prediction errors are  also developed.
Given the parametric model,  the proposed methodology 
automatically produces estimates of the small area 
parameters and their MSPEs, without   requiring
explicit analytical  expressions for the MSPE.\\[.1in]
\noindent
{\it Some key words:} Best predictor; Bootstrap; Kernel; Mean squared
 prediction error.

\section{\small\sc Introduction}

Small area estimation (SAE)  is an important statistical 
research area due to its growing demand from public 
and private agencies. The variance  of a small area estimator 
based on the direct small area  sample is unduly  large and hence,
 there is a need for  constructing model 
based estimators with low mean squared prediction error
(MSPE).  
 A good account
 of small area estimation research is available in a
 recent book by J.N.K. Rao (Rao, 2003). Although,
in theory,  it is possible to use  such a model based 
 approach, in practice a statistician  often faces  
some  challenging problems in implementing it  due to the 
fact that  for each model, estimators must be 
derived and their properties studied.
Indeed, a small deviation from the standard model assumptions 
may require a considerable amount of analytical work 
and need special expertise. For example, Prasad and
 Rao (1990)  (hereafter referred to as PR) derived small area estimation formulas 
assuming normality of both the sampling distribution 
and the population distribution (for two-level small 
area models, discussed later) and with the moment 
based estimators of model parameters. After about a
 decade, Datta and Lahiri (2000) extended this 
approach when the model parameters are estimated by
the  maximum likelihood approach. Recent  works 
of Jiang, Lahiri and Wan (2002)  and Lahiri and Maiti 
(2003) (hereafter referred to as JLW and 
LM, respectively) allow a more general framework, but both works
require  the knowledge of the exact functional forms
of the MSPE, which are known   only  in few  simple cases. 
 However, a general  solution to finding the best 
estimator of the small 
area parameters or of its functions, and estimation 
of the associated MSPE  are  not 
available. A second problem with the existing approaches
(except for the LM method)  to estimating the MSPE is 
that these methods do not
 always produce non-negative estimates. { Though the 
linearization technique of 
PR  produces non-negative estimates under normality, 
the jackknife method may produce negative MSPE estimates (Bell, 2002).}
 Consequently,
 there is a great demand for a general estimation 
system where the user can only specify the distributions
 and then valid estimates of the small area parameters
 and their MSPEs  can be obtained without much
 of analytical efforts. 

In  this paper, we consider a general two level aggregate data model
and  develop a unified system  for prediction  of 
 small area parameters and estimation of the associated 
MSPE.  Here we  extend  the ``perturbation'' or ``tilting" method 
of LM   and  
construct  a nonnegative   estimator of the MSPE   that 
achieves  second order accuracy for bias 
correction without  requiring explicit 
analytical derivation of the MSPE function.
The key idea is to combine the LM approach 
with the parametric version of the  bootstrap method of Efron (1979)
so that accurate numerical approximations to 
various intermediate population quantities can be 
generated numerically. We  show that under 
some  regularity conditions, the proposed MSPE estimator 
   attains  second order accuracy
%The most important feature of the proposed methodology,
%which sets it apart from the existing methods of
%bias corrected MSPE estimation,  is that it produces
%nonnegative, second order accurate  estimators
%of the MSPE  
for  a wide  range of parametric distributions 
and for a general class  of model parameter estimates
and their nonlinear functions, without  requiring the user to 
derive the formulas on a case by case basis.

The rest of the paper is organized as follows. In Section 2,
we consider the general two level aggregate data model that is commonly
used in the context of small area estimation. In Section 3,
we describe the (estimated)  best predictor of  
 functions of the small area parameters. In Section 4, we 
briefly describe the  existing  approaches to MSPE estimation
 and  also give a description of the proposed method.
Theoretical properties of the proposed method are given in 
Section 5. Results from a  simulation study 
and some concluding remarks  are presented in Sections 6 and 7, respectively.
Proofs are given in  the Appendix.

\section{\small \sc Generalized Mixed  Models for Small Area Estimation}

Consider the general two level aggregate data model
%\begin{equation}
\be
%\left\{\begin{array}{ll}
y_i|\theta_i  \stackrel{ind} {\sim} F_1(\cdot; \theta_i,R_i),\quad  
\theta_i  \stackrel{ind} \sim  F_2(\cdot;{\ux}_i,\ulambda, G_i), ~ i=1,\cdots,m, 
\ee
 where, $R_i$ and $G_i$ are known functions of a vector of $p$-parameters 
$\upsi=(\psi_1,\ldots,\psi_q)$,
say, $(R_i,G_i) = g_i(\upsi)$. Thus,
the model is determined by the parameter vector 
 $\udelta\equiv (\ulambda^T,\upsi^T)^T$,
 a $(p+q)\times 1$ vector of constants. Usually, 
$y_i$'s are direct survey estimators with sampling
 variance $R_i$, $\theta_i$'s are small area
 parameters and $\ux_i$ a set of covariates 
available at the estimation stage. Aggregate   
and generalized
linear mixed effects models are special cases of (2.1).

Consider the Fay-Herriot (1979) type small area model
\be
y_i=\theta_i+e_i, \quad
\theta_i=\ux_i^T\ulambda+v_i
\ee
where $e_i$'s are independent $N(0,s_i)$ with known 
$s_i$, $v_i$'s are iid
$N(0,\sigma_v^2)$ and $e_i$ and $v_i$'s are
 independent. Furthermore, 
$\ux_i$ is a known $p\times 1$ vector of 
co-variates, $\ulambda$ is the
vector of regression coefficients; $y_i$ 
is the direct survey estimator of
$\theta_i$. Note that (2.2) can be written
 as
$
y_i=\ux_i^T\ulambda+v_i+e_i
$
which is a special case of a linear mixed model where both $F_1$ and 
$F_2$ are normal cdf.

%\begin{equation}
%\uy=\uX\ulambda+\uZ\uv+\ue
%\end{equation}
%where $\uy=(y_1,\ldots,y_m)^T, \uv=(v_1,\ldots,v_m)^T\sim 
%N(0,\sigma_v^2I_m)$
%and $\ue= (e_1,\ldots,e_m)^T\sim N(0,R)$, with
% $R= diag(s_1,\ldots,s_m)$ { and $\uZ$ is a known
%diagonal  matrix}. 
%Consequently (2.2) and (2.3) are special cases of (2.1) 
%where  both $F_1$ and $F_2$ are normal.

Next consider the mixed logistic model, where conditional on 
small area parameter $p_i,$ the direct estimator  $y_i$ is 
binomial $(n_i,p_i), i=1,\ldots,m$; here $n_i$ is the number 
of sampled units in the
$i$-th small area. Then, consider the model 
\be
\theta_i\equiv logit(p_i)=\ux_i^T\ulambda+v_i,
\ee
 where the  $v_i$'s
are iid $N(0,\sigma_v^2)$. 
In this case $F_1$ is binomial and $F_2$ is normal in the logit scale.
This is a special case 
of generalized linear
mixed model.
% that has received considerable attention in recent
% years, e.g.
%Breslow and Clayton (1993), and Lee and Nelder (1996).
% Ghosh {\it et al.} (1998)
%used a more general version of  model (2.4) for disease mapping. 

Our objective is to make inference about a function of the  small area
 parameter $\theta_i$
\be
\beta_i=h(\theta_i), i=1,\ldots,m,
\ee
where $h$ is a suitable function chosen by the user. 
%For example, in 
%the case of the 
%logit-normal model, one may be interested in predicting the mean
%response $p_i=logit^{-1}(\ux_i^T\ulambda+v_i)$. For 
%the normal-normal model,  often a log-transformation is applied.
 For example, the ``Small Area Income and Poverty Estimation''
(SAIPE)  project
 of the US Census Bureau uses the  log value of the direct
 estimates for estimating poverty at the county level and thus an inverse 
transformation needed for the parameter of interest.
We  would like to emphasize that, at the second level of modeling,
 the structure always need not  be  of the form
$
h(\theta_i)=\ux_i^T\ulambda+v_i
$. 
In fact, we can also use nonlinear modeling,  such as 
$
h(\theta_i)=\kappa(\ux_i;\ulambda,v_i). 
$ where $\kappa$ is a nonlinear function.
%For example,  in the case of a  nonlinear mixed effect model, we may have 
%$\kappa(\ux_i;\ulambda,v_i)=g(\ux_i^T\ulambda)+v_i$,
%where $g(\cdot)$ is some nonlinear function.

\section{\small \sc Development of the Best and Empirical Best Predictors}
\setcounter{equation}{0}
As an estimator of the small area parameter, we will
 take the best 
predictor (BP) as defined below. 
%The best predictors are 
%equivalent 
%to best linear unbiased predictors (BLUPs)  in
% some special  cases, such 
%as in the  normal-normal case. 
 Let $\beta_i=h(\theta_i)$ 
be the 
 parameter 
of interest. We define the BP  and the  
{\it empirical best predictor} (EBP) of  $\beta_i$, respectively, by  
\begin{eqnarray}
%eqn(3.1)
\tilde{\beta}_i
&=& E_{\udelta}\{h(\theta_i)|\uy\},\\
\hat{\beta}_i&=&E_{\hat{\udelta}}\{h(\theta_i)|\uy\},
\quad  1,\ldots, m,
\end{eqnarray}
where
$\hat{\udelta}$ is an estimator of  $\udelta$. 
 For example, in the Fay-Herriot model
(2.2), the BP of $h(\theta_i) =\theta_i$ 
takes the form 
$
\tilde{\beta}_i=\ux_i^T\ulambda+\frac{\sigma_v^2}{\tau_i}
(y_i-\ux_i^T\ulambda)$,
where $\tau_i=\sigma_v^2+s_i$. For a general $h(\cdot)$, however, 
a closed form simple expression for the BP/EBP  and their  MSPE 
may not be  available. Consequently, the PR-type SAE methodology 
based on Taylor's expansions may not be readily  applicable.

Next, we derive some  useful general
 formulas for the EBP of (3.2). Note that by  the 
independence of  $y_i$'s, 
 the  conditional distribution of $\theta_i$ given $y_1,
\ldots,y_n$ depends only on $y_i$ (and $\delta$). Hence, 
$
\tilde{\beta}_i =
 E_{\udelta}\{h(\theta_i)|\uy\}
= \int h(t)F_{\theta_i|y_i}(dt;\udelta) 
\equiv \xi_i(y_i;\udelta)
$ 
{say},
where $F_{\theta_i|y_i}(\cdot;\udelta)$ denotes the 
conditional distribution of $\theta_i$ given $y_i$.
The EBP is  given by
\be  
%%%%%%%%%%%%%%%%%%%%%%%%%%%%%%%%%%%%%%%%%%%%%%%%%%%%%%%%%%%%%%%%%%%%%%eqn (3.3) -r.2
\hat{\beta}_i =
\xi_i(y_i;\hat{\udelta}).
\ee
First  consider the case where the marginal distribution 
of $\theta_i$ has a probability density function (pdf)
$f_2(\cdot;\ux_i,\ulambda, G_i)$  (with respect to the 
Lebesgue measure) and the conditional distribution
$F_1(\cdot;\theta_i,R_i)$ of $y_i$ given $\theta_i$ has
a generalized  density 
$f_1(\cdot;\theta_i,R_i)$ 
(i.e., the Radon-Nikodym derivative
 with respect to a $\sigma$-finite measure). For example,
$f_1$ can itself be a pdf or a probability mass function
(pmf) for a discrete probability distribution. In this case, 
the EBP 
%function $\xi_i(y_i;\udelta)$ defining 
%the best predictor $\tilde{\beta}_i$ in (3.3) 
is given by
\be
%%%%%%%%%%%%%%%%%%%%%%%%%%%%%%%%%%%%%%%%%%%%%%%%%%%%%%%%%%%%%%%%%%%%%%eqn (3.4) -r.2
\hat{\beta}_i = \xi_i(y_i;\hat{\udelta})
=\frac{\int h(t) p_i(y_i,t;\hat{\udelta}) dt}
{\int  p_i(y_i,t;\hat{\udelta}) dt},
\ee
where $p_i(y,t;\udelta)= 
f_1(y;t,R_i) f_2(t;\ux_i,\ulambda, G_i)$.
 Next   consider the case where  the marginal distribution 
of $\theta_i$ is discrete and has a pmf 
$f_2(\cdot;\ux_i,\ulambda, G_i)$  and $F_1(\cdot;\theta_i,R_i)$
has a generalized  density  $f_1(\cdot;\theta_i,R_i)$ as above.
Here 
the EBP
%
%$\xi_i(y_i;\udelta)$  of (3.3) 
%
 is given by
\be
%%%%%%%%%%%%%%%%%%%%%%%%%%%%%%%%%%%%%%%%%%%%%%%%%%%%%%%%%%%%%%%%%%%%eqn (3.5)  -r.2
\hat{\beta}_i
=
\xi_i(y_i;\hat{\udelta})
=\frac{\sum_{t} h(t) p_i(y_i,t;\hat{\udelta})}
{\sum_{t}   p_i(y_i,t;\hat{\udelta})},
\ee
where  $p_i(y,t;\udelta)$ is  as before 
%= 
%f_1(y;t,R_i) f_2(t;\ux_i,\ulambda, G_i)$
 and where the sum
in (3.5) runs over all $t$ in the support of $\theta_i$.
In many applications, formulas (3.4) and (3.5) can be
implemented  using numerical methods,
e.g., numerical integration, MCMC, importance sampling, etc.
  For example,
for the logit-normal model with 
the canonical link, 
  $\xi_i(y_i;{\udelta})
=\big[\int \alpha_i(t)(y_i+1)\{1+\alpha_i(t)\}^{-n_i-1}\phi(t)dt\big]/\big[
\int \alpha_i(t) y_i \{1+
\alpha_i(t)\}^{-n_i}\phi(t)dt\big]
$, 
where $\alpha_i(z)=\exp(\ux_i^T\ulambda + \sigma_{v}z)$ and $\phi$ is the 
N(0,1) pdf (e.g., see, McCulloch and Searle (2001, pp 273) and JLW). In this case, the EBP  can be
easily evaluated by generating N(0,1) 
variates and using the Monte-Carlo method.

\noindent
{\bf Remark  1:} ({\it Parameter estimation}). In general,  the 
maximal likelihood estimates (MLE's)  do not have any closed form 
expressions. Except for  the conjugate and linear 
link models, maximization of the marginal likelihood involves integration with respect 
to the distribution function $F_2$. There is no unique way of  evaluating  
this integral. Using  advanced  techniques  such as EM based 
MLE, Markov Chain Monte Carlo (MCMC) based MLE, etc.,  the MLE's can be
 obtained for a  large  class of distributions. An excellent account of guidelines 
for the general mixed linear models can be obtained in Chapter 10 of McCulloch 
and Searle (2001). We mention that the SAE  methodology 
developed here is equally applicable for other type of parameter estimators 
such as those based on method of moments or estimating equation approaches,
provided  
they are ${m}^{1/2}$ consistent. 

\noindent
{\bf Remark  2:} ~For  situations where a direct implementation 
of (3.4) or  (3.5) is  difficult, we now describe some approximations 
to the EBP using the bootstrap method of Efron (1979) and 
the nonparametric functional estimation
methodology.
%
%
%From a user point of view, it would be  very useful
% to have a unified solution,  where one does not need to
% derive  the exact expressions for the EBP and its MSPE 
%analytically on a case by case basis. 
%This motivates us 
%to develop based algorithm, which work for a general
% class of distribution $f_1(\cdot)$ and $f_2(\cdot)$.
%
%
Note that $\hat{\beta}_i$ is the conditional expected
 value of a function of $\theta_i$ for 
fixed $\delta$  evaluated at $\delta=\hat{\delta}$.
%
%given the observations $\uy$.
 This suggests that 
under mild regularity conditions,
we may  approximate  $\hat{\beta}_i$
to any desired level of accuracy by using 
standard regression 
function estimation methods, such as 
Nadaraya-Watson estimators, local polynomial estimators, etc.
%
%without having to find out an exact expression for 
% $\hat{\beta}_i$. 
%
Let  $\{y_i^{*j},
\theta_i^{*j}\}_{j=1}^J$ be generated  values using  model (2.1), 
but with   $\udelta=
\hat{\udelta}$.
When the distributions of $\theta_i$ and $y_i$ are   continuous, 
 we propose a Nadaraya-Watson 
 approximation to $\hat{\beta}_i$, given by
\begin{equation}
%%%%%%%%%%%%%%%%%%%%%%%%%%%%%%%%%%%%%%%%%%%%%%%%%%%%%%%%%%%%%%%%%%%%%%% eqn (3.6) -r.2
\hat{\beta}_i^*
=\frac{\sum_{j=1}^Jk_{(\cdot)}(y_i^{*j}-y_i)h(\theta_i^{*j})}
{\sum_{j=1}^Jk_{(\cdot)}(y_i^{*j}-y_i)},
\end{equation}
where $k_{(\cdot)}$ is a symmetric kernel  function chosen suitably. 
 There are many choices of $k_{(\cdot)}$, such as 
a Gaussian kernel $k_{(\cdot )}(x)=\frac{1}{b}k(x/b)$
where  $b$ is the bandwidth and $k(x)=\phi(x)$, the standard
normal density function.
On the other hand,  when the marginal distribution of $y_i$ is
 discrete,  
%and $f_1$ is a pmf, 
we propose
\begin{equation}
%%%%%%%%%%%%%%%%%%%%%%%%%%%%%%%%%%%%%%%%%%%%%%%%%%%%%%%%%%%%%%%%%%%%%%% eqn (3.7) -r.2
\hat{\beta}_i^*
=\frac{\sum_{j=1}^Jh(\theta_i^{*j})I_{[y_i^{*j}=y_i]}}
{\sum_{j=1}^JI_{[y_i^{*j}=y_i]}},
\end{equation}
where $I_{[\cdot]}$ denotes the indicator function. 
Results on Nadaraya-Watson estimators  of regression 
functions imply  (cf. H\"ardle (1991)) that
%given $\uy$, 
\begin{equation}
%%%%%%%%%%%%%%%%%%%%%%%%%%%%%%%%%%%%%%%%%%%%%%%%%%%%%%%%%%%%%%%%%%%%%%% eqn (3.8) -r.2
|\hat{\beta}_i - \hat{\beta}_i^*|^2 
     = O \Big( (Jb)^{-1} + b^{-2}\Big)
\qmq{in  probability,}
\end{equation}
% 
%%%%%%%%%%%%%%%%%%%%%%%%%%%%%%%
%%%%%%% in conditional probability, to be precise
%
 as $J\raw \infty$ and $b\raw 0$
 in such a way  that $Jb\raw \infty$.
The bound in  (3.8) is available uniformly
over $i=1,\ldots,m$, provided there exists a constant $C\in(0,\infty)$
such that 
$E_{\udelta}|\xi_{i}^{''}(Y_i;\ut)| + 
E_{\udelta}|g_i^{''}(Y_i;\ut)| < C
$ 
for all
$i=1,\ldots,m$ and for all $\ut\in {\cal N}$, a neighborhood
of the true value of the unknown parameter $\udelta$. Here,
 $\xi_{i}^{''}(y;\ut)
=\frac{\partial^2}{\partial y^2}\xi_{i}(y;\ut)$,  
$g_{i}^{''}(y;\ut)
=\frac{\partial^2}{\partial y^2}
g_{i}(y;\ut)$, 
and $g_{i}(y;\ut)$ is the marginal density of $Y_i$.
%(given by the denominator of (3.5)).
For the discrete case, a direct 
computation shows that
%given $\uy$, 
\begin{equation}
%%%%%%%%%%%%%%%%%%%%%%%%%%%%%%%%%%%%%%%%%%%%%%%%%%%%%%%%%%%%%%%%%%%%%%% eqn (3.9) -r.2
|\hat{\beta}_i - \hat{\beta}_i^*|^2 
     = O \Big( J^{-1}\Big)
\qmq{in  probability,}
\end{equation}
as $J\raw \infty$. This  bound is also available  uniformly
 in $i$, provided 
$E_{\udelta} |h(\theta_i)|^2 
%+ E_{\udelta}|g_{i}(Y_i;\ut)|^2
 <C
$
 for all $i$ and for all $\ut\in {\cal N}$,
where $C\in (0,\infty)$  is a constant, and  ${\cal N}$
is as above.
%$g_{i}(y;\ut)$ is the marginal 
%density  (pmf) of  $Y_i$, given by the denominator 
%of (3.6).

% 
%%%%%%%%%%%%%%%%%%%%%%%%%%%%%%%
%%%%%%%%% in conditional probability, to be precise
%
%
%%
%The following result shows that a sufficiently close 
%approximation to $\hat{\beta}_i$ can be generated using the 
%above formulae.\\[.2in]
%
%\noindent
%{\bf Theorem 1:} (a) Suppose that $\xi_i(y;\bd)$ is twice
%continuously differentiable in $y$ with 
%$
%|\frac{\partial^1 \xi_i(y;\bd)}{\partial y^2}| \leq G_0(y;\bd)$
%for some nonnegative function $G_0$ with $E G_0(y_i;\hud)=O(1)$,
% $f_2$ is continuous and positive on its support, and condition (C.2)
%of Section 5 holds.
%Then, there exists a constant $C_1\in(0,\infty)$ such that 
%$$
%E|\hat{\beta}_i - \hat{\beta}_i^*|^2 \leq 
%    C_1 \Big[ (Jb)^{-1} + b^{-2}\Big],
%$$
%
%as $J\raw \infty$ and $b\raw 0$. 
%
%(%%b) Suppose that $f_2$ is discrete and condition (C.2)
%of Section 5 holds..  Then 
% there exists a constant $C_2\in(0,\infty)$ such that 
%$$
%E|\hat{\beta}_i - \hat{\beta}_i^*|^2 \leq 
%     C_2 J^{-1},
%$$
%as $J\raw \infty$.\\
%
Thus,  for both the discrete and the continuous data,
  the accuracy of the approximation
$\hat{\beta}_i^*$ to $\hat{\beta}_i$  increases
with larger values of $J$.  For the continuous case,
 we need to specify a choice of
the bandwidth $b$. For  kernels arising from symmetric probability
densities, the optimal choice of  $b$ is of the order  $J^{-1/5}$. 
We take  the bandwidth $b$ of this optimal order, e.g.,  $b=J^{-1/5}$, and 
%Since we can generate an arbitrarily large number of samples 
%from the estimated distribution, we can guarantee
attain a desired level of accuracy by choosing $J$ suitably large.
% The details of implementation of the 
% procedure will be discussed in
% the implementation section. 
Finite sample accuracy of the approximations (3.6) and (3.7) 
are typically very good. See Table 1 in Section 6 below which
 reports    the relative biases and MSPE'S of (3.6) and (3.7) for
the  normal-normal 
and the logit-normal examples.  

%In the next section, we describe a general 
%methodology for assessing
% the accuracy of the EBPs. 

\section{\small \sc Mean Squared Prediction Error and its Estimation}
\setcounter{equation}{0}
\subsection{\small\it Background}
As a measure of accuracy of the EBP $\hat{\beta}_i$,
we shall consider  the
Mean Squared Prediction Error(MSPE) of $\hat{\beta}_i$, ~ 
$
MSPE(\hat{\beta}_i)=E_{\udelta}(\hat{\beta}_i-\beta_i)^2 \equiv M_{i}(\delta)
$.
It is easy to show  that 
\begin{eqnarray}
%%%%%%%%%%%%%%%%%%%%%%%%%%%%%%%%%%%%%%%%%%%%%%%%%%%%%%%%%%%%%%%%%%%%%%% eqn (4.1) -r.2
M_i(\delta)
= E_{\udelta}(\tilde{\beta}_i-\beta_i)^2+E_{\udelta}(\hat{\beta}_i
-\tilde{\beta}_i)^2 
\equiv M_{1i}(\udelta)+M_{2i}(\udelta), \quad\mbox{say}.
\end{eqnarray}
%
%The above decomposition is standard and each 
% term on the
%ight side of (4.1) has a clear interpretation. 
The first term 
$M_{1i}(\delta)$ is the mean
squared error of the (ideal) best predictor $\tilde{\beta}_i$ 
while the second term  $M_{2i}(\udelta)$ accounts for  the extra 
variability due to the  estimation of $\udelta$.  
Typically, 
%Furthermore,  under some 
%mild regularity conditions, 
\be
%%%%%%%%%%%%%%%%%%%%%%%%%%%%%%%%%%%%%%%%%%%%%%%%%%%%%%%%%%%%%%%%%%%%%%% eqn (4.2) -r.2
M_{1i}(\udelta)=O(1)\qmq{and}
M_{2i}(\udelta)=O(m^{-1})
\qmq{as}\mti.
\ee
It is tempting  to plug in $\hat{\udelta}$ in (4.2) 
and get a simple 
MSPE estimate as
\begin{equation}
%%%%%%%%%%%%%%%%%%%%%%%%%%%%%%%%%%%%%%%%%%%%%%%%%%%%%%%%%%%%%%%%%%%%%%% eqn (4.3) -r.2
mspe_{\mbox{\sc sim}}(\hat{\beta}_i)
=M_{1i}(\hat{\udelta})+M_{2i}(\hat{\udelta}) .
\end{equation}

%
%{\bf Theorem:}
%
% \begin{eqnarray*}
%E M_{1i}(\hat{\delta})&=&M_{1i}(\delta)+0(\frac1m)\\
%E M_{2i}(\hat{\delta})&=&M_{2i}(\delta)+o(\frac 1m)
%\end{eqnarray*}
%
%I don't know how these results associate with kernel width and
% no of iteration on J.
%
%It follows from the above theorem that 

However, this approach has  two drawbacks.
First, explicit expressions for the functions $M_{1i}(\udelta)$
and $M_{2i}(\udelta)$  are not always available.
In the very special case of  the 
 normal-normal Fay-Herriot model, an  expression for $M_{1i}(\udelta)$
 and an approximation for  $M_{2i}(\delta)$ are available
for $h(\theta_i)= \theta_i$, $i=1,\ldots,m$. 
%{\it In fact $M_{2i}(.)$ 
%was further partitioned into two pieces (see, 
%Rao, 2003 equ. 6.2.31)}. 
Even for this model,   expressions are 
not available  for a {\it nonlinear}  function of $\theta_i$ and one has 
to derive those. For example, Slud and Maiti (2006) 
derived the expressions for MSPE estimates under normal 
set up when $h$ is an exponential function. 
%
%The expressions under 
%the  logit normal model
%(even for $h(\theta_i)= \theta_i$, $i=1,\ldots,m$) 
%are very difficult to derive. 
%

The second problem with the above  approach is a little more subtle. 
To describe it, note that 
typically, the estimator $\hat{\udelta}$  has  bias and variance   of
 order $O(m^{-1})$, which  propagate  through the simple  MSPE estimator,
 leading to 
$
E \{M_{1i}(\hat{\udelta})\}=M_{1i}(\udelta)+O(m^{-1})\qmq{and}
E \{M_{2i}(\hat{\udelta})\}=M_{2i}(\udelta)+o(m^{-1})
$
as $\mti$. (Here and in the following, we  often  drop the 
subscript $\udelta$ to ease notation). 
Thus, 
$
E\{ M_{1i}(\hat{\udelta})- M_{1i}(\udelta)\}
$, 
the bias of 
the simple  estimator  of $\moi(\bd)$,  
is of the order $O(m^{-1})$  %. 
%Consequently, for the naive estimator 
%of MSPE($\hat{\beta}_i$), 
%the bias of 
%$M_{1i}(\hat{\udelta})$ 
 which  masks the contribution
of $ M_{2i}(\cdot)$ to the MSPE of $\hat{\beta}_i$
(cf. (4.2)).

In view of the second problem,
 in the SAE literature,  it is customary 
 to require that 
the bias of a ``good''
estimator of MSPE$(\hat{\beta}_i)$  be of {\it smaller order than} 
$O(m^{-1})$.  
Traditionally, the  bias of the naive estimator
$M_{1i}(\hat{\udelta})$
  is reduced   by explicit bias correction, either by using 
a Taylor's expansion of the function $M_{1i}(\cdot)$ (cf. PR)
%Prasad and Rao (1990))
or using the Jackknife method (cf. JLW). 
%
%Jiang, Lahiri, and Wan (2002))
% (hereafter referred to as PR and JLW, respectively). 
%Because of the explicit bias adjustments, 
% the resulting MSPE estimators  can take negative values.  
%{ Though under normality, the PR method of 
% bias correction always produces nonnegative MSPE estimates,
% the jackknife method may produce negative estimates (Bell, 2002).} 
%%%%%%%%
%
% ----> Move to the comparision section! 
%
%{\it However, for applying PR  method, one needs to know the closed form
% expressions  of the functions  $M_{1i}(\cdot)$ 
%(and  its derivatives) and 
% $M_{2i}(\cdot)$ . 
%%%%%%%%%%%%%%%%
Other related work include 
Pfeffermann and Tiller (2005) and  Pfeffermann and Glickman (2004).
The first paper   approximated 
$M_{2i}(.)$ and the bias of $M_{1i}(\hat{\udelta})$ under a 
state space model based on parametric bootstrap,  assuming normality of the errors.
The second used  a bias corrected estimator of $M_{1i}(.)$ and a parametric bootstrap 
estimator of $M_{2i}(.)$, for the Fay-Herriot model. Pfeffermann and Glickman 
also developed  a `nonparametric' bootstrap method that did  not require 
generating samples from a distribution. Nonetheless, normality was still 
assumed implicitly.
% since they used  the facorization (4.1).
 In a recent work, LM 
  proposed a new 
approach to  bias correction that attains second order  
accuracy  and at the same time, produces a  nonnegative estimator of the MSPE. 
%  For the 
%jackknife method, JLW (2002) mentioned at least twice that they 
%required a closed form expression for $M_{1i}(.)$ and 
%approximated the $M_{2i}(.)$ using jackknife variance 
%estimation. Rao (2003) laid out an approach of implementing 
%jackknife method for binary data, but is not clear for 
%general distribution.} 
%
%
Here, 
we extend the 
LM approach    to the case of estimating the 
MSPE of a general  function of $\theta_i$ with 
second order accuracy under a general two-level parametric  model,
 even when  exact expressions
for the functions $M_{1i}(\cdot)$ and 
 $M_{2i}(\cdot)$ are not available. 
For completeness, we now briefly describe the LM method.  
Suppose that for $i=1,\cdots,m$, 
\vspace*{-.1in}
\be
%%%%%%%%%%%%%%%%%%%%%%%%%%%%%%%%%%%%%%%%%%%%%%%%%%%%%%%%%%%%%%%%%%%%%%% eqn (4.4) -r.2
 \sum_{j=1}^k|M_{1i}^{(j)}(\udelta)| >\epsilon_0, 
\vspace*{-.1in}
\ee
for some $\epsilon_0 >0$, 
where  for a  smooth  function $f: \Re^k\rightarrow \Re$, 
$f^{(j)}$, $f^{(j,r)}$ and $f^{(j,r,s)}$  denote the first,  the second 
and the third 
order partial derivatives with respect to the $j$-th 
co-ordinate, the $(j,r)$-th co-ordinates,
and the  $(j,r,s)$-th co-ordinates, respectively, 
$j,r,s=1,\cdots,k$, where $k$ is the number of model parameters.  
Condition (4.4) says that  $M_{1i}^{(j)}(\udelta)\neq 0 $ for some  $j$. 
For notational simplicity,  we suppose that $M_{1i}^{(1)}(\udelta)\neq 0$. Then,
 the {\it preliminary perturbed estimator} of $\udelta$ for 
the $i$-th small area is defined as  
$
\bar{\udelta}_i  =  \hat{\udelta}-\hat{B}_i
\{M_{1i}^{(1)}(\hat{\udelta})
\}^{-1}\ue_1
$
where
$
 \hat{B}_i \equiv 
\sum_{j=1}^kM_{1i}^{(j)}
(\hat{\udelta})\hat{b}(j)
+\frac 12 \sum_{j=1}^k\sum_{r=1}^kM_{1i}^{(j,r)}
(\hat{\udelta})\hat{V}(j,r)
$, 
  with    $\hat{\ub}=(\hat{b}(1),\ldots,\hat{b}(k)$ and 
 $\hat{V}=(( \hat{V}(j,r)))_{k\times k}$ 
respectively denoting  some suitable estimators (e.g., bootstrap estimators)
of the bias and the variance of
$\hat{\udelta}$,
%%%%%%%%%%%%%%%%%%%%%%%%%%%%
%. Thus, the estimator
%$
%\bar{\udelta}_i
%$
%is obtained from the initial estimator $\hat{\udelta}$
%by adding a correction factor to the first component of
%$\hat{\udelta}$ only. Note that if instead of 
%$M_{1i}^{(1)}(\udelta)$, a different partial derivative 
%$M_{1i}^{(r)}(\udelta)$ were nonzero, then we would
%define the  preliminary tilted  estimator 
%$\bar{\udelta}_i$ by replacing the factor 
%$\left\{M_{1i}^{(1)}(\hat{\udelta})
%\right\}^{-1}\ue_1$ in (4.5) with 
%$\left\{M_{1i}^{(r)}(\hat{\udelta})
%\right\}^{-1}\ue_r$, 
%where the vector
and  $\ue_r \in \Re^k$
has $1$ in the $r$-th position and zeros elsewhere,
$1\leq r\leq k$. 
%%%%%%%%%%%%%%%%%%%%%%%%%%%%%%%%%%%%%%%%
%
%Next, let $\Delta$ denote the set of possible values of the parameter 
%$\udelta$ under  model (2.1). Then 
%$i=1,\cdots, m$. Thus, if the preliminary estimator $\bar{\udelta}_i$ 
%takes values inside the parameter space $\Delta$ and the 
%value of the partial derivative $M_{1i}^{(1)}(\hat{\udelta})$
%at $\hat{\udelta}$ is not too small, the perturbed 
%estimator of $\udelta$ is given by $\bar{\udelta}_i$ itself. 
%However, in the event that either $\bar{\udelta}_i$ falls outside 
%$\Delta$ or $M_{1i}^{(1)}(\hat{\udelta})$ becomes too small, 
%we replace it with the original estimator 
%$\hat{\udelta}$. Small values of $M_{1i}^{(1)}(\hat{\udelta})$
%make the estimator $\bar{\udelta}_i$ unstable and hence,
%these are 
%truncated below. LM show  that under 
%appropriate regularity conditions, the probability of 
%getting a preliminary estimator $\bar{\udelta}_i$ 
%outside $\Delta$ or that of getting a value of
%$M_{1i}^{(1)}(\hat{\udelta})$ below the threshold
%$(1+\log m)^{-2}$ tends to zero rapidly 
% as $m\rightarrow \infty$, uniformly in $i$. As a 
%consequence, the perturbed 
%estimator $\check{\udelta}_i$ coincides with the 
%preliminary perturbed  estimator $\bar{\udelta}_i$ 
%with high probability. 
%%%%%%%%%%%%%%%%%%%%%%%%%%%%%%%%%%%%%%%
%
The LM  estimator of 
the MSPE is now defined as 
\be
%%%%%%%%%%%%%%%%%%%%%%%%%%%%%%%%%%%%%%%%%%%%%%%%%%%%%%%%%%%%%%%%%%%%%%% eqn (4.5) -r.2
mspe_{\mbox{\sc lm}}(\hat{\beta}_i) 
= M_{1i}(\check{\udelta}_i)+M_{2i}(\hat{\udelta}), 
i=1,\cdots, m, 
\ee
where $\check{\udelta}_i$ is the {\it perturbed estimator} 
 of $\udelta$ for the $i$-th small area, defined by 
\be
%%%%%%%%%%%%%%%%%%%%%%%%%%%%%%%%%%%%%%%%%%%%%%%%%%%%%%%%%%%%%%%%%%%%%%% eqn (4.6) -r.2
\check{\udelta}_i = \left\{\ba{ll}
        \bar{\udelta}_i & \mbox{ if } \bar{\udelta}_i \in \Delta 
\mbox{ and } |M_{1i}^{(1)}(\hat{\udelta})|^{-1} \le (1+\log m)^2 \\
    \hat{\udelta} & \mbox{ otherwise,} \ea \right.
\ee
and $\Delta$ is  the set of possible values of the parameter 
$\udelta$ under  model (2.1).
Note that by  construction, the MSPE estimator is 
always {nonnegative}. Further, LM  show that under 
some regularity conditions, the  bias  of  
the estimator $mspe_{\mbox{\sc lm}}(\hat{\beta}_i) 
$
is of the order $o(m^{-1})$.
% and a standard error that is of the order
% $O(m^{-\frac12})$. 
% Therefore, the proposed estimator attains 
%the same level of accuracy as the previously proposed 
%estimators by  PR and JLW, while at 
%the same time, it guarantees non-negativity.  

\noindent
{\bf Remark 3:} When more than one partial derivatives 
$\moi^{(j)}(\udelta)$ are non-zero, one may use 
perturbations along all such directions.
Thus, an  alternative  MSPE estimator is
 given by
\be
%%%%%%%%%%%%%%%%%%%%%%%%%%%%%%%%%%%%%%%%%%%%%%%%%%%%%%%%%%%%%%%%%%%%%%% eqn (4.7) -r.2
mspe_{\mbox{\sc  LM:alt}}(\hat{\beta}_i) = \moi(\check{\udelta}^{\dagger}_i)
    + \mtoi(\hat{\udelta}_i), ~i=1,\cdots, m,
\ee
where 
$$
%\be
\check{\udelta}^{\dagger}_i = \left\{\ba{ll}
        {\udelta}^{\dagger}_i & \mbox{ if } {\udelta}^{\dagger}_i \in \Delta 
\mbox{ and } |{\cal J}|^{-1}
 \sum_{j\in {\cal J}}|M_{1i}^{(j)}(\hat{\udelta})|^{-1} \le (1+\log m)^2 \\
    \hat{\udelta} & \mbox{ otherwise,} \ea \right.,
%\ee
$$
$
\udelta^{\dagger}_i = \hat{\udelta}- \sum_{j\in {\cal J}}
\big[\hat{B}_i/\moi^{(j)}(\hat{\udelta})\big]
 \ue_j/|{\cal J}|
$,   ${\cal J}=\{j: 1\leq j\leq k, \moi^{(j)}(\udelta) \neq 0\}$
and for any set $A$, let $|A|$ denotes its size.  
The arguments developed in LM readily imply that the new MSPE estimator
is also second order correct, under the same set of regularity conditions
 as in  LM. By combining all $|{\cal J}|$ directions, the new
estimator attains a better finite sample stability.
% (as  the chances of $\check{\udelta}^{\dagger}_i\neq 
%{\udelta}^{\dagger}_i $ become   smaller). 

%%%%%%%%%%%%%%%
%
%Note that in defining the modified preliminary estimator $\udelta^{\dagger}_i$,
%we have distributed the total `bias correcting factor' 
%$ \hat{B}_i$ (cf. (4.5)) {\it equally} among the $|{\cal J}|$ directions. 
%Arguably,  a better   alternative is to  use a {\it  weighted}  average %of $\hat{B}_i$
%along  the $|{\cal J}|$ non-zero components,
%with weights proportional to the magnitude of $\moi^{(j)}(\udelta)$. Under  this
%scheme, 
%  components with smaller values of $[\moi^{(j)}(\hat{\udelta})]^{-1}$
%get higher weights.  The optimal allocation of the weights is attained
%if we again concentrate all the weights along a  component $j_0$  that minimizes  
% $[\moi^{(j)}(\hat{\udelta})]^{-1}$.  This yields
%a solution  similar to the original proposal of LM of using a single
%co-ordinate to make the bias corection, but with the important 
%difference that the selected co-ordinate $j_0$ is now data-dependent.
%From a computational point of view, neither of the  modifications 
%require significantly higher efforts (than using a preselected $j$) 
%since  all the partial derivatives   $[\moi^{(j)}(\hat{\udelta})]^{-1}$ 
%must be   computed to find   $\hat{B}_i$, in any case. 
%
%%%%%%%%%%%%%%%%%%

\subsection{\small \it Nonnegative estimation of the MSPE  when expressions for 
$\moi$ and $\mtoi$ are Unavailable}
As discussed  earlier, except for very few  
standard models, exact or closed form expressions  for  the terms 
 $\moi(\cdot)$ and $\mtoi(\cdot)$ are not available.
Here we employ the Bootstrap method of Efron (1979) to
develop  an  approximated version  of the estimator
  mspe$_{\rm LM}$ that 
is nonnegative, second order accurate, and that can
be computed without additional analytical work on
 the part of the user. 
%%%%%%%%%%%%%
To that end, first we define 
a bootstrap based approximation to the function $M_{1i}(\cdot)$
at a given value $\udelta_0$.  Let $(y_i^{*l}, 
\theta_i^{*l}), ~l=1,\ldots,N_0$ be iid random vectors
generated using  model (2.1) with $\udelta = \udelta_0$. Then the 
bootstrap approximation to $M_{1i}(\udelta_0)$ is 
given by
\be
%%%%%%%%%%%%%%%%%%%%%%%%%%%%%%%%%%%%%%%%%%%%%%%%%%%%%%%%%%%%%%%%%%%%%%%%%% eqn (4.8)
M_{1i}^*(\udelta_0) = 
\frac{1}{N_0}\sum_{l=1}^{N_0} 
		\Big \{\xi_i(y_i^{*l};\udelta_0)
 - h(\theta_i^{*l})\Big\}^2.
\ee
%where the dependence of $M_{1i}^*(\udelta_0)$ on  $N_0$ 
%is suppressed for notational convenience.
Next we use $ M_{1i}^*(\cdot)$  to 
construct  estimators  of   
 the partial derivatives of 
the function $M_{1i}(\cdot)$.
%%%%%%%%%%%%
% at $\udelta$ using the 
%approximating function $M_{1i}^{*}(\cdot)$. 
%%%%%%%%%%%%%
To motivate the 
construction,
 consider a smooth  function  
$f:\Re \rightarrow  \Re$. Then, for any 
$a\in\Re$,
\begin{eqnarray*}
f(a+\epsilon) -f(a-\epsilon)
= \{f(a+\epsilon) -f(a)\}
	-
		\{f(a-\epsilon) -f(a)\}
= 2\epsilon f'(a) + o(\epsilon)
\end{eqnarray*}
{as} $\epsilon\raw 0$,
where $f'(a)$ denotes the derivative of $f(\cdot)$ at $a$.
Hence the scaled difference $(2\epsilon)^{-1} 
\big\{
f(a+\epsilon) -f(a-\epsilon)\big\}$
gives an approximation to $f'(a)$
for small values of $\epsilon$. We now employ this fact to define 
suitable approximations to the first order partial derivatives 
of $M_{1i}(\cdot)$ at $\hat{\udelta}$. Let $\{\z\}$
be a sequence of positive real numbers converging to
zero. 
 Let
\be
 %%%%%%%%%%%%%%%%%%%%%%%%%%%%%%%%%%%%%%%%%%%%%%%%%%%%%%%%%%%%%%%%%%%%%%%%%% eqn (4.9)
M_{1i}^{(j)*}(\hat{\udelta}) 
= \frac{1}{2\z}\Big\{ M_{1i}^*(\hat{\udelta} + \z\ue_j)
	- M_{1i}^*(\hat{\udelta} - \z\ue_j)
	\Big\},
\ee
$j=1,\ldots,k$. 
 Using a similar reasoning, we also define 
approximations to the second order partial derivatives as 
% of$M_{1i}(\cdot)$ at $\hat{\udelta}$ by
\begin{eqnarray}
%%%%%%%%%%%%%%%%%%%%%%%%%%%%%%%%%%%%%%%%%%%%%%%%%%%%%%%%%%%%%%%%%%%% eqn (4.10)
M_{1i}^{(j,j)*}(\hat{\udelta}) 
&=& \frac{1}{[\z]^{2}}\Big\{ M_{1i}^*(\hat{\udelta} + \z\ue_j)
	+ M_{1i}^*(\hat{\udelta} - \z\ue_j)
	- 2M_{1i}^*(\hat{\udelta})
	\Big\},  ~~1\leq j \leq k,\\
%%%%%%%%%%%%%%%%%%%%%%%%%%%%%%%%%%%%%%%%%%%%%%%%%%%%%%%%%%%%%%%%%%%% eqn (4.11)
M_{1i}^{(j,r)*}(\hat{\udelta}) 
&=& \frac{1}{2{\z}^{2}}\Big[
	\Big\{
		 M_{1i}^*(\hat{\udelta} + \z\ue_{j,r})
	+ M_{1i}^*(\hat{\udelta} - \z\ue_{j,r})
	- 2M_{1i}^*(\hat{\udelta})
		\Big\}	\nonumber\\
	&&\hspace{.6in}
		- \z^{2}\Big\{
		M_{1i}^{(j,j)*}(\hat{\udelta}) 
		+ M_{1i}^{(r,r)*}(\hat{\udelta}) 
				\Big\}
	\Big], ~~1\leq j\neq r\leq k,
\end{eqnarray}
where $\ue_{j,r} = \ue_j + \ue_r$.  Theorem 1 in Section 5 shows that
under some regularity conditions,
$
\max_{1\leq i\leq m} E| M_{1i}^{(j)*}(\hud) -
 M_{1i}^{(j)} (\hud)| = O\big( \z
	+(\z)^{-1}N_0^{-{\eta}/({1+\eta})}\big).
$
and 
$
\max_{1\leq i\leq m} E| M_{1i}^{(j,r)*}(\hud) -
 M_{1i}^{(j,r)} (\hud)| = O\big( \z
	+(\z)^{-2}N_0^{-{\eta}/({1+\eta})}\big).
$
for all $1\leq j,r\leq k$, for some  $\eta\in(0,1]$.
 Thus, the proposed estimators of the 
partial derivatives provide accurate approximations 
for  suitable choices of $\z$ and $N_0$. 
%
%\raw \infty$ in such a way that $N_0^\frac{\eta}{1+\eta}} 
%[\z]^2\raw \infty$. 
%

Next for $l=1,\ldots, N_0$, let $(y_1^{*l},\ldots, y_m^{*l})$
be iid  random vectors having  joint  distribution 
(2.1) with $\udelta =\hat{\udelta}$ and let  
 $\udelta^{*l}$  
denote the bootstrap version  of $\hat{\udelta}$, obtained by 
replacing $(y_1,\ldots,y_m)$ with $(y_1^{*l},\ldots, y_m^{*l})$.
Define  the bootstrap estimators of the bias and the variance of 
$\hat{\udelta}$  by
\be
%%%%%%%%%%%%%%%%%%%%%%%%%%%%%%%%%%%%%%%%%%%%%%%%%%%%%%%%%%%%%%%%%%%% eqn (4.12)
\ub^* = \frac{1}{N_0}\sum_{l=1}^{N_0} \udelta^{*l} -\hat{\udelta}
~~\mbox{and}~~  V^* =  \Big\{\frac{1}{N_0}\sum_{l=1}^{N_0} \udelta^{*l}(\udelta^{*l})^T\Big\}
		- \Big(\frac{1}{N_0}\sum_{l=1}^{N_0} \udelta^{*l}\Big)
		\Big(\frac{1}{N_0}\sum_{l=1}^{N_0} 
				\udelta^{*l}\Big)^T,
\ee
respectively. 
Theorem 2 in Section 5 below  gives conditions for the consistency of
 $\ub^*$ and $V^*$.
%
%The proposed estimator of the MSPE uses a version of the 
%tilted estimator of $\udelta$ to accomplish  an implicit bias correction
%as in LM. Indeed,
%
With this, we  now  define   the {\it bootstrap based   preliminary perturbed 
estimator} $\bar{\udelta}_i^{*}$ as 
$$
\bar{\udelta}_i^*  =  \hat{\udelta}- B_i^*
\big\{M_{1i}^{(s)*}(\hat{\udelta})
\big\}^{-1}\ue_s,\qmq{provided} |M_{1i}^{(s)}(\udelta)| \neq 0
\qmq{for some} 
s\equiv s_i\in\{1,\ldots,k\},
$$
  where  
$ B_i^* =
\sum_{j=1}^k M_{1i}^{(j)*}
(\hat{\udelta}){b}^*(j)
+2^{-1} \sum_{j=1}^k\sum_{r=1}^kM_{1i}^{(j,r)*}
(\hat{\udelta}){V}^*(j,r)
$,
%%%%%%%%%
%and $\epsilon_1\in (0,1)$, where 
%$\epsilon_1$ does not depend on $i,m$. 
%%%%%%%%%%
%Here, and in the following,
 $x(j)$ denotes the $j$th component of a vector  $\ux$ and 
$B(j,r)$ denotes the $(j,r)$th element of a matrix  $B$.
The {\it bootstrap based perturbed estimator} of $\udelta$ for the 
$i$th small area is now defined as
\be
%%%%%%%%%%%%%%%%%%%%%%%%%%%%%%%%%%%%%%%%%%%%%%%%%%%%%%%%%%%%%%%%%%%% eqn (4.13)
\check{\delta}_i^* = \left\{\ba{ll}
        \bar{\udelta}_i^* & \mbox{ if } \bar{\udelta}_i^* \in \Delta 
\mbox{ and } |M_{1i}^{(s)*}(\hat{\udelta})|^{-1} \le (1+\log m)^2 \\
    \hat{\udelta} & \mbox{ otherwise } \ea \right.
\ee
 and the bias corrected estimator of $M_{1i}(\udelta)$
is given by
$
M_{1i}^*(\check{\delta}_i^*).
$,
$i=1,\cdots, m$.

Next we define the bootstrap estimator of $M_{2i}(\udelta)$.
Note that 
$M_{2i}(\udelta) = E_{\udelta}(\hat{\beta}_i - 
\tilde{\beta}_i)^2 = E_{\udelta}\big\{\xi_i(y_i,\hat{\udelta})
- \xi_i(y_i,{\udelta})\big\}^2
$. 
 Let 
 $\udelta^{*l}, ~l=1,\ldots, N_0$ 
denote iid bootstrap replicates  of $\hat{\udelta}$ as above (cf. (4.12)). 
% denote iid random vectors 
% having   distribution  (2.1) with $\udelta 
%= \hat{\udelta}$.
%
% Let $\udelta^{*l}$ denote the bootstrap
%version of $\hat{\udelta}$ based on $(y_1^{*l},
%\ldots,y_m^{*l})$,  for $l=1,\ldots, N_0$. 
Then, the parametric   bootstrap 
estimator of $M_{2i}(\udelta)$ is now defined as
\be 
M_{2i}^{*}(\hat{\udelta}) 
	= N_0^{-1}\sum_{l=1}^{N_0} \big\{ 
\xi_i(y_i^{*l},{\udelta}^{*l})
- \xi_i(y_i^{*l},\hat{\udelta})\big\}^2.
\ee
Pefferemann and Tiller (2005), Pfeffermann and Glickman (2004) and 
Butar and Lahiri (2003) also proposed similar parametric 
bootstrap estimates of $M_{2i}(.)$ for normal errors.
% when the error distributions are normal. 

The proposed {\it bias corrected estimator} of the MSPE
$M_{i}(\udelta)$ is defined as
\be
mspe_{\mbox{{\sc new}}}(\hat{\beta}_i) = M_{1i}^* (\check{\delta}_i^*)
+
M_{2i}^{*}(\hat{\udelta}),
\ee
$i=1,\cdots, m$.
%Note that by the construction, the proposed MSPE estimator is 
%always nonnegative. Further, to find  the proposed estimator,
%
In the next section, we show that under 
some regularity conditions, the proposed estimator 
 has a bias that is of the  order $o(m^{-1})$.
% and a standard error that is of the order
% $O(m^{-\frac12})$. 
As a result, the proposed estimator attains 
the same level of asymptotic bias  accuracy as the previously proposed 
MSPE  estimators. Furthermore, as (4.15) does not require 
 explicit expressions 
for  the functions $M_{1i}$ and $M_{2i}$,
the proposed MSPE estimation methodology 
 can be applied to complex or nonstandard  models where none of the existing methods 
are easily applicable.  

\section{\small \sc Theoretical Properties}
\setcounter{section}{5}
\setcounter{equation}{0}
For investigating the theoretical properties
of the proposed method,  we shall suppose that the 
random variables $(y_i,\theta_i): i=1,\ldots,m$
and the various bootstrap variables 
$(y_i^{*l},\theta_i^{*l})$'s are defined on a common 
probability space. We write $P_{\ux}$ and $E_{\ux}$ to denote the 
probability and the expectation under a given parameter 
value $\ux\in\De$. For notational simplicity, we set 
$P_{\bd} = P$ and $E_{\bd} = E$ where $\udelta$ is 
the true value of the parameter.  Define the functions
$\ua(\cdot)$ and $\Sigma(\cdot)$ by 
$\ub(\bd)   = \ua(\udelta)/m$
and $V(\udelta)  =\Sigma(\udelta)/m$,
where 
$\ub(\bd) \equiv E_{\udelta}(\hud) -\udelta$
and $V(\udelta) = var_{\udelta}(\hud)$. Note that
$\ua$, $\ub$, $\Sigma$ and $V$ depend on  
  $m$. 
Unless otherwise specified, limits in
the order symbols below are taken as $\mti$.
 Also, let $\ecd$ denote the conditional
expectation of the bootstrap variables given $\hat{\udelta}$.
Proofs of the main results are given in the Appendix.
\\[.1in]
\noindent
{\bf Conditions}
\begin{itemize}
\item[(C.1)] ~(i)  ~$\udelta$, the true value of the parameter, is an 
interior point of $\Delta$. \\[.1in]
(ii) $M_{1i}$ is three times  continuously differentiable on 
$\Delta$ and there exists a constant $C_1\in (0,\infty)$ 
such that for all $\ux \in \Delta, j,r,s=1,\cdots,k$ and $i=1,\cdots,m, m\geq 1$,
$|M_{1i}^{(j)}(\ux)|+|M_{1i}^{(j,r)}(\ux)| +
|M_{1i}^{(j,r,s)}(\ux)|
<C_1 
$.\\
(iii)  $M_{2i}$ is differentiable on $\Delta$ and 
 there exist constants $C_2, C_3, \epsilon_0  \in (0, \infty)$ 
and $\gamma \in (0,1]$ and a function 
$G_1:\Re^k\rightarrow [0,\infty)$ with $E G_1(\hat{\udelta})=O(1)$
such that for all $j=1,\cdots,k; i=1,\cdots,m, m\geq 1$, ~~
%\begin{eqnarray*}
$|M_{2i}^{(j)}(\udelta)| \leq C_2 m^{-1}$,
$
|M_{2i}(\ux)| \leq  m^{-1}G_1(\ux) \qmq{for all} \ux \in \Delta
$ {and}
$m|M_{2i}^{(j)}(\ux)-M_{2i}^{(j)}(\udelta)|\leq  C_3\|\ux-\udelta\|^{\gamma}
$
{for all} $ \ux \in \N$, 
%\end{eqnarray*}
where $\N\equiv \{\|\ux-\udelta\|\leq \epsilon_0\}$. 
\item[(C.2)]
There exist constants $\eta \in (0,1]$ 
and $C_4\equiv C_4(\eta) \in (0, \infty)$ such that 
$
%\sup_{\udelta\in\Delta} (Add?)
E|h(\theta_i)|^{2+2\eta} < C_4$
{for all}  $i=1,\cdots,m, m\geq 1
$.
\item[(C.3)] 
(i) ~Let $\rho_m(\ux; a) = E_{\ux}\|\hud -\ux\|^a$, $\ux\in\Delta, a\in (0,
\infty)$. Suppose that there exists a constant $\eta \in (0,1]$ 
such that 
$
E\rho_m(\hud; 2+2\eta) = O(1)
$.\\
(ii) The sequences of  functions $\{\ua\}\equiv \{\ua_m\}$ and 
$\{\uSi\}\equiv\{\uSi_m\}$
are (component-wise) equicontinuous at $\bd$.\\
(iii) There exists a continuous function $G_2:\Re^k \raw [0,\infty)$
such that $\|\ua_m(\ux)\| + \|\uSi_m\| \leq G_2(\ux)$
for all $\ux\in\uDelta$ and $EG_2(\hud)^2 =O(1)$.
% [FIX!]
\item[(C.4)] 
There exists a  constant $\eta \in (0,1]$  such that
$E\Big\{\sqrt{m}|\hud -\bd\|\Big\}^{2+2\eta} =O(1)$.
\end{itemize}

We now briefly comment on the regularity conditions.
Condition (C.1) requires the functions $M_{1i}$ and $M_{2i}$
to be smooth, which typically holds under suitable smoothness conditions
on the parametric model (2.1). As mentioned
earlier, in most applications the function $M_{1i}$ is of the 
order $O(1)$ while $M_{2i}$ is of the 
order $O(m^{-1})$ as $\mti$.  Condition (C.1)  requires that the 
partial derivatives of these functions also have the 
same orders. 
Conditions (C.2), (C.3)(i), and (C.4) are
moment conditions depending  on  $\eta$, whose values will
be specified in the statements of the results below.
These  are used to prove `closeness' of
various parametric bootstrap estimates to 
their conditional expectations.  Condition (C.3)(ii) and (iii)
are exclusively used to establish consistency of the bootstrap estimators
of the bias and the variance estimators of $\hud$.

The first result proves consistency of the partial derivative 
estimates.\\[.1in]
\noindent
\noindent
{\bf Theorem 1:} Let  Conditions (C.1)(ii) and (C.2) hold and let  
$N_0$ be as in   (4.8). Then
\begin{eqnarray}
\max_{1\leq j\leq k}
\max_{1\leq i\leq m} E\Big | M_{1i}^{(j)*}(\hud) -
 M_{1i}^{(j)} (\hud)\Big| &=& O\Big( \z
	+[\z]^{-1}N_0^{-\frac{\eta}{1+\eta}}\Big),\\
\max_{1\leq j,r\leq k}
\max_{1\leq i\leq m} E\Big | M_{1i}^{(j,r)*}(\hud) -
 M_{1i}^{(j,r)} (\hud)\Big| &=& O\Big( \z
	+[\z]^{-2}N_0^{-\frac{\eta}{1+\eta}}\Big).
\end{eqnarray}

Note that the right sides of (5.1) and (5.2) go to zero
for any  
$\z\raw 0$, $N_0\raw \infty$ such that $\z^2N_0^{{\eta}/({1+\eta})}
\raw \infty$.  Here $\z$ acts as a smoothing parameter 
that controls the bias parts of the proposed estimates.
For a  smaller value of $\z$, a larger value of  the resample size $N_0$
has to  be chosen accordingly to attain a desired accuracy level.
Also, note  that the value of $N_0$
 required for estimating the second order partial derivatives 
must  grow at a faster rate than the case of the first order partial derivatives
to attain the same level of accuracy.

The next result considers accuracy of the bootstrap
bias and variance estimators of $\hud$. \\[.1in]
\noindent
{\bf Theorem 2:} Let  Condition (C.3) hold and let 
$N_0$ be as in  (4.8). Then, 
$$E\|\ub^* -\hub\|^2 = O(N_0^{-1})\qmq{and} 
E\|V^* -\hV\| = O(N_0^{-\frac{\eta}{1+\eta}}).
$$
 Under the conditions of Theorem 2, the bootstrap bias estimator
is ${N_0}^{1/2}$-consistent. The variance estimator can also attain 
the same rate, provided $\eta=1$. Note that unlike Theorem 1,
the estimators of the bias and the variance matrix of $\hud$
do not involve a smoothing parameter like $\z$.

The next result shows that 
under suitable conditions, the proposed estimator of the
MSPE($\hat{\beta}_i$)  second order bias corrected. \\[.1in]
\noindent
{\bf Theorem 3:} Suppose that conditions (C.1)-(C.4) hold 
and that $\eta=1$ in both (C.2) and (C.3). Suppose that 
for each $i=1,\ldots,m$, there exists $s=s(i)\in\{1,\ldots,k\}$
such that 
\be
|\moi^{(s)}(\bd)|>C_0
\ee
for all $m\geq 1$, where $C_0\in(0,\infty)$.
 Let $\z=m^{-5/4}$ 
and $N_0\geq m^{a}$ for some $a>9/2$.
%%%%%%%%old  $a>7/2$ - fails on the bound on R_{2i}.
%\begin{enumerate}
%\item[(i)]
% Then
%\be
%\max_{1\leq i\leq m} E \Big|mspe_{\mbox{{\sc new}}}(\hat{\beta}_i)
%	 - mspe_{\mbox{{\sc lm}}}(\hat{\beta}_i)\Big|
%=o(m^{-1}).
%\ee
%\item[(ii)]
%Further, 
Then  the proposed  mspe estimator is second order bias 
accurate, i.e., 
\be
 \max_{1\leq i\leq m} \Big| 
	E \Big[mspe_{\mbox{{\sc new}}}(\hat{\beta}_i)\Big] 
- M_i(\udelta)\Big|
=o(m^{-1}).
\ee
%\end{enumerate}

 Theorem 3 shows that the proposed MSPE estimator
 achieves the same 
 second order bias accuracy as the earlier methods proposed in the literature.
Thus, under the given regularity  conditions,
 the additional randomness induced by several 
resampling steps  has a negligible effect on the bias of 
the new estimator.
Since it also does not require the knowledge
of the functions $\moi(\cdot)$, $\mtoi(\cdot)$, of their   the 
partial derivatives,  and of the bias and variance of the 
estimator $\hud$, the proposed method can be applied to
any model of the form (2.1), where the other methods
are not readily applicable. The price paid for this 
omnibus solution 
%to the second order bias corrected  MSPE estimation
%problem
 is that it is  computationally  intensive.
% nature of the proposed method. 

%In the next section, we explore  the 
% performance of the proposed method and compare it 
%with the existing approaches, whenever available.

\section{\small\sc Practical Implementation and Numerical  Findings}
\subsection{\small \it Finite sample considerations}
\setcounter{equation}{0}
In this section, we provide some guidelines for implementing the 
proposed MSPE estimation methodology in finite sample applications.
Supposing, for the time being, that an expression for the BP is 
known, computation of different parts of  the estimator 
$mspe_{\mbox{\sc new}}$ involves 
generating (parametric) bootstrap samples from the joint distribution
of $(y_i,\theta_i)$ for $i=1,\ldots,m$ (cf. (2.1)) at various values of 
the parameter $\udelta$.  For the bootstrap bias and variance estimators
$\hub$ and $\hV$ and the term $\mtoi^*(\hud)$, 
we suggest using a  resample size (drawn   from (2.1) with $\bd =\hud$) 
 in the    100s (e.g., in the range 500-1000). This 
is known to be adequate for Monte-Carlo evaluation of bootstrap
estimators of  variance-type functionals (cf. Efron and Tibshirani (1993)).
Next consider  numerical evaluation of the first term of 
$mspe_{\mbox{\sc new}}$, i.e., of   $\moi^*(\check{\udelta}_i^*)$.
This requires us to approximate   the partial derivatives of
$\moi(\cdot)$  which, in turn,   involve the smoothing parameter $z_m$. 
For all computations done in this section,  we 
 set $z_m=m^{-5/4}$ as in Theorem 3, although other choices of 
$z_m \ll m^{-1}$ may be used. For the numerical approximation
of the partial derivatives,   the resample sizes must  be larger
in order to compensate for the effect of the smoothing - the smaller the
choice of the  smoothing parameter $z_m$, the larger the choice of $N_0$ 
will have to be. For $m$ of moderate size (e.g., $m\in (10, 80)$)
and $z_m$ as above, we have found 
resamples of size $N_0$ in the range 2000-10,000  adequate for computing 
the first order partial derivatives $\moi^{*(j)}$ and resamples of size
$N_0 \approx  10,000$   for the second order ones $\moi^{*(j,r)}$. 
Finally, in the case that an exact expression for the  EBP is not available
and it is approximated numerically using (3.6) or (3.7), 
the resample size $J$  may be chosen in the  100s
(e.g., 300-1000) in the discrete case while  it must be of a higher 
order (e.g., 1000+) in the continuous case.
Approximations  given by the above choices of the resample sizes
are generally very good.  In the next section, we 
report the results of a simulation study and the associated computing time
for     three  specific examples 
where we follow the finite sample guidelines given above. 
For an illustration,  
Table 1 below gives  the resulting  approximations 
 for the EBP both in the discrete and the 
continuous cases which appear to be in good agreement with 
the true values.

\subsection{\small \it Simulation results}
In this section, we check the performance of  the  
 MSPE estimators (4.5) and  (4.15) for Models I-III described below,
 and compare them  with the Datta-Lahiri (2000) (hereafter, referred to as DL)
 version of the PR method  and the jackknife method of JLW, 
 as described in Rao (2003). 
DL extended the PR method 
when the  model parameters are estimated using 
MLEs. 
 We used MLEs of the model parameters for Models I and II,
and used estimating equations  for Model III.  Normal kernel was 
used for the kernel based EBP's.
 We use the following notations 
for different methods of MSPE estimation: ~ 
JK for  jackknife, LM1 for (4.5) and LM2 for (4.15). 
%We consider the following situations.
\\[.1in]
\noindent
%\begin{description}
%\item[{\bf Model I: Normal-Normal.}]
%\vspace*{.1in}
{\bf Model I: Normal-Normal.}
~This is a continuous data model, where  both $F_1$ and $F_2$ are normal; 
The model structure is   specified by  (2.2)
%$$
%y_i= \theta_i + e_i\\
%\theta_i = \ux_i^T\ulambda + v_i
%$$
%where $e_i$'s are independent $N(0,s_i)$ with known $s_i$
%and $v_i$ are iid $N(0,\sigma^2_v)$, $e_i$ and $v_i$s are independent,
%and  
 with $\ulambda =0$. In this setting, 
all four methods of bias correction  are applicable.
Although in this case a closed form expression  for the 
BP exists, to gain some insight into the performance of the 
suggested  approximations, we use (3.6) to find the BP 
for the LM2 method.
For the other  three  methods, the available closed form 
expressions are  used.   
We choose  $F_2$ to be  normal with mean 0 and variance unity,
and  $F_1$ to be  
 normal with mean 0 and variance $s_i, i=1,\cdots, m.$ with  m=15. 
The 15 areas  are divided into three groups of five, 
with equal numbers of areas and equal values of 
$s_i$. The three different values of $s_i$ used are 
$(.7, .5, .3)$. The set-up is 
similar to the one  considered by Datta {\it et al.} (2005). \\[.1in]
%We denote this 
%continuous data model  as $M_1$. 
%%%%%%%%%%%%%%5
%\item[{\bf Model II: Binomial-`Logit-Normal'.}]
%\vspace*{.1in}
\noindent
{\bf Model II: Binomial-`Logit-Normal'.}
~This is a binary data model where 
we suppose that  $F_1$ is binomial and $F_2$ is logit-normal. 
In particular, the logit of the success probability 
 of $F_1$ is normally distributed 
with mean zero and variance unity (cf. (2.4) with $\ulambda=0$). In this 
setting, only JK and LM2 methods of MSPE estimation are applicable. 
The binomial  population has 8 areas, of respective sizes 
$n_i$=36, 20, 19, 16, 17, 11, 5  and 6, based on the number of 
patients receiving  a particular treatment from different 
clinics (Booth and Hobert, 1998). To generate the $i$th binomial 
population, we first 
generate the success probability
\be
%%%%%%%%%%%%%%%%%%%%%%%%%%%%%%%%%%%%%%%%%%%%%%%%%%%%%%%%%%%%%  eqn (6.1) -r2 
p_i=\frac{\exp(\mu+v_i)}{1+\exp(\mu+v_i)}
\ee
where $v_i$ is a standard normal variate, $i=1,\cdots,8$ and $\mu=0$.
 In this case  the BP is not available in a closed form. We first find the maximum likelihood 
estimates of the model parameters using Slud (2000). Then the BP is calculated 
using Gauss-Hermite quadrature with 15 points for the JK method and (3.7) 
for the LM2 method. \\[.1in]
%We denote this binary data model as $M_2$.
%%%%%%%%%%%%%%%%%%%%%
%\item[{\bf Model III: Normal-Lognormal.}]
%\vspace*{.1in}
\noindent
{\bf Model III: Normal-Lognormal.} 
~This is a continuous data non-conjugate model, where  $F_1$ is normal and $F_2$ is 
 lognormal. You and Rao (2002) considered this model for estimating the Canadian census
 undercoverage and called this as `unmatched sampling and linking model'. 
%However,their  approach was fully Bayesian and 
Here, neither the PR/DL nor the JK methods are   applicable in a straightforward way.
 We took $m=15$ and generated  $\theta_i$'s $(i=1,\cdots,m)$ from a lognormal 
distribution. We took two covariates, besides the intercept, one was 
generated from $N(0,.5)$ and the other was generated 
from Uniform $(.5,1)$. We set $\udelta=(\ulambda^T,\sigma_v^2)^T=(0,0.5,-1.5,0.5)^T$.
 Then given $\theta_i$'s, $y_i$'s were generated as in  Model I. Instead of using 
ML estimate, we used unbiased estimating equation approach 
for estimating the model parameters (cf. Ghosh and  Maiti, 2004). 
%since it is easy to obtain the 
%moments of the marginal distribution of $y$ using the properties of normal and 
%lognormal distributions and conditioning arguments. 
%
Since the BP does not have any closed form expression, we used the kernel based 
estimator (3.6) for estimating the BP and consequently, of the four, here 
LM2 is the only method available for estimating the MSPE. 
Also, note that in this case, one can obtain the {\it perturbed estimator} of
 $\udelta$ either by (4.15) or by the (estimated version of the) method 
described in Remark 3. Both  methods gave very similar   results.   
 The  MSPE estimator  in Remark 3 (with estimated  partial derivatives, etc.)
gave slightly low CV than the estimator in  (4.15); see Table 2. 
%\end{description}

In implementing LM2, we used 1000 bootstrap samples for finding the bias and 
variances 
%of the model parameter 
estimates and 10000 bootstrap samples for all 
other approximations. All  simulation results were based on
 R=1000 replication.  The approximate computation time for each 
model is at most  48 hours on  a 
UNIX workstation equipped with 4000MHz 64-bit CPU and FORTRAN 77 compiler. 
In any real application,  user needs to run the code only once, meaning minimal computational 
time (less than 3 minutes)  with  data sets of a  similar size. 

To study the performance of the EBP $\hat{\theta}_i$ of the  small area 
parameter $\theta_i$, we use the following two empirical measures.
\bea
\mbox{{\it Absolute relative bias}}\hspace{.2in}
T_1&=&\frac{1}{R}\sum_{r=1}^R|\frac{\hat{\theta}_i^{(r)}
-\theta_i^{(r)}}{\theta_i^{(r)}}|.
\\
\mbox{
{\it Empirical MSPE}}
\hspace{.2in}
T_2 &=& \frac{1}{R}\sum_{r=1}^R(\hat{\theta}_i^{(r)}
-\theta_i^{(r)})^2.
\eea
The body of all the tables gives averages over all the small areas where 
the ``average" is measured in terms of the median (given in the first column 
for each model) or the mean (in the second column).

\newpage
\begin{center}
 Table 1.
{\it  Absolute relative bias ($T_1$) and empirical MSPE 
($T_2$) for the EBP. Results using the 
   kernel based approximations (3.6) and (3.7)   are reported within the parentheses.}
% under the
% results using standard methods.

\begin{tabular}{ccccccc}\\ 
       &\multicolumn{2}{c}{Model I}&\multicolumn{2}{c}{Model II} 
& \multicolumn{2}{c}{Model III} \\ 
Measures & Median & Mean & Median & Mean & Median & Mean\\ 

$T_1$    & 2.318  & 4.171 & 0.223 & 0.243 & --- & --- \\ 
         &(2.156) &(4.122)&(0.224)&(0.243)& (0.926) & (0.923) \\
$T_2$    & 0.376  & 0.366 & 0.0107& 0.0131 & --- & ---\\
         &(0.378) &(0.373)&(0.0107)&(0.0131)& (0.269) & (0.292)\\ 
\end{tabular}
\end{center}
\vspace*{2ex}

There is a good agreement between the actual values
 and the approximations for the EBP given in equations
(3.6) and (3.7). For the binary data, this agreement is 
particularly remarkable. This is because for the same value 
of the resample size $J$, the approximation in the discrete case
is more accurate (having a faster  rate of convergence). 
 In the case of the binary data, the ``actual'' values
are found by numerical integration.  
 The simulation result 
shows that 
both the numerical integration based approximation   
and the ``kernel''  method based  approximation (3.6)
 behave similarly. However, kernel method seems 
more automated  than  numerical  integration 
as it does not require additional programming for a different
continuous data  model.
% 
%from the user's perspective.

Table 2 reports the  following empirical measures of relative bias 
and coefficient of variation,  
 quantifying  the performances of different MSPE estimation methods:
\bea
\mbox{{\it Relative bias}}
\hspace{.2in}
T_3 &=& [E\{\hat{MSPE}(\hat{\theta}_i)\}-T_2]/{T_2}\\
\mbox{{\it Coefficient of variation}}
\hspace{.2in}
T_4&=&\left[E\{\hat{MSPE}(\hat{\theta}_i)
-T_2\}^2\right]^{\frac 12}/{T_2}.
\eea
Here $E\{\hat{MSPE}(\hat{\theta}_i)\}$ 
and $E\{\hat{MSPE}(\hat{\theta}_i)-T_2\}^2$ are  estimated 
empirically by   averaging  the replicates of 
 $\hat{MSPE}(\hat{\theta}_i)$ and  $\{\hat{MSPE}(\hat{\theta}_i)-T_2\}^2$,
respectively.

\newpage

\begin{center}
 Table 2.
{\it  Relative biases ($T_3$) and coefficient of variations ($T_4$)
 for the bias corrected estimators of the MSPE.
Entries within parentheses represent LM1 and LM2 estimates based on Remark 3 modification.
}

\begin{tabular}{cccccccc}\\ 
       &   &\multicolumn{2}{c}{Model I}&\multicolumn{2}{c}{Model II} 
& \multicolumn{2}{c}{Model III} \\ 
Method & Measures & Median & Mean & Median & Mean & Median & Mean \\ 

PR/DL  & $T_3$    & -0.016  & -0.004 &  ---  & ---  & --- &---\\
       & $T_4$    & 0.159   &  0.150    &  ---  & ---  & --- & ---\\
JK     & $T_3$    & 0.068   &  0.095    & -0.088 & -0.026& --- & --- \\
       & $T_4$    & 0.504   &  0.635    &  0.686 &  0.758& --- & --- \\
LM1    & $T_3$    & -0.015  & -0.018    &  ---  &  --- & --- & ---\\
       &          &(-0.000) & (0.050)   &  ---  &  --- & --- & --- \\
       & $T_4$    & 0.158   &  0.151    &  ---  &  --- & --- & ---\\
       &          &(0.153)  & (0.149)   &  ---  &  --- & --- & --- \\
LM2    & $T_3$    & -0.013  & -0.028    & -0.108 & -0.083&  0.116&  0.041  \\
       &          &(-0.019) &(-0.024)   &(-0.087)&(-0.083)&(0.115)&(0.044) \\
       & $T_4$    & 0.229   &  0.224    &  0.172 &  0.164& 0.319  & 0.368 \\
% \hspace{0in}0.319 (0.310)& \hspace{0in} 0.368 (0.298) \\
       &          &(0.225)  & (0.218)   &(0.170) &(0.156)&(0.310) &(0.298)\\
       
\end{tabular}
\end{center}

For Model I, all the methods perform well in terms of minimizing relative bias. 
However, in terms of the coefficient of variation, there is a difference in the 
performance of the four methods. The PR/DL and LM1 
methods turn out to be the best, followed by the  LM2 method.
 The small increase in the  variation of the LM2 method over the
LM1 method  is expected, as the randomness in  the various
approximation  steps in its construction adds  to the total variability
of the bias corrected MSPE estimator. However, the highest 
variation for  this model is observed for the JK method, where 
 the  variation more than {\it double} compared to the LM2 method
and it is more than {\it three} times compared to the LM1 and PR/DL 
methods.% This behavior is sustained also 

As mentioned earlier, for Model II,  only the LM2 and the JK methods 
are applicable. In this case,  the  LM2 tends to have higher relative bias.
However, in terms of the coefficient of variation, which gives the 
{\it combined} effects of the bias and the variance of the MSPE estimators,
the LM2 method again beats the JK method by a relative magnitude  of 
 300\% to 400\% or more. To gain further insight into the 
bias properties of the two methods, 
%Given the number of small area is 8, this may not be a  serious 
%concern. To confirm this, 
we repeated the  simulation study with $m=16$  areas (instead of the
$m= 8$  areas considered earlier) under Model II. For this higher  
value of $m$, we  found that the relative bias for  the LM2 method
dropped to -.038 and  -.024 for the median and the mean, respectively. 
The eight additional small area sizes were 37, 32, 19, 17, 12, 10, 9 and 7.
In comparison, the relative bias for  the JK method under $m=16$
were -.025 and -.026 for the median and mean respectively. The 
coefficient of variations for the two methods continued to show 
 a  similar pattern as in the $m=8$ case.
Thus for both  models, the estimators 
 produced by the JK method  have  inferior performance in terms of
the  coefficient of variation.
% The bias corrected estimator given by DL method (PR type) and LM1 have
% the  lowest relative bias and the coefficinet of variation. The LM2 method 
% continue to enjoy low relative bias with an increasing coefficient of 
% variation but significantly less than the jackknife method.

For Model III, the PR/DL method is not applicable and the existing 
literature does not show how to apply the JK method. This is a 
somewhat unusual set up of simulation within the existing
SAE  literature. It may be interesting to know that, if some one naively 
used  $M_{1}(\hat{\udelta})$ with  formula (4.8), the median relative bias
would be   -.227 and  the mean,  -.260. This indicates severe under-estimation 
which is expected. In comparison, 
 LM2 produces satisfactory results for both the relative bias and 
the coefficient of variation.

\section{\small \sc Discussion}%%%%%%%%%%%%%%%%%%%%%%%%%%%%%%%%%%%%%%%% Section 7
In this paper, we consider a new method of bias correction 
for the ``simple'' estimator of the MSPE of a 
possibly nonlinear function  of the small area means
 $h(\theta_i),i=1,\cdots,m$.  The proposed method may be contrasted 
with the existing methods, which require explicit  analytical expressions
for bias correction. The  popular 
linearization method of bias correction proposed by PR  
can not be easily extended to  nonlinear $h$ and non-normal models. 
Further the PR approach is sensitive to the method of estimating 
model parameters in the sense that additional 
analytical work may be 
needed for each new estimation method. 
In the cases where  exact analytical expressions are available, 
 the simulation results indicate that the PR and the LM methods 
(are comparable and) have the best overall performance (in terms of MSEs)
while 
the proposed   method  (LM2)  fares reasonably well against these.
In particular,  LM is preferable to  LM2 in such situations.  
As for comparison with the JK method in this case, 
the  LM2 method performs  much better than the 
JK method in finite samples in terms of the co-efficient of variation. 

In the more complicated examples, where exact analytical expressions
for the MSPE are not available, the LM and the PR  methods are 
not applicable, but the LM2 method and the JK method
(with some suitable  adaptation)  are. In  this case,  the LM2 method  
seems to have a superior performance compared to   the JK 
 method in terms of overall accuracy. 
%We also mention that 
Further, because of the  inherent limitations of the Jackknife method for estimating 
the variance of a non-smooth estimator of the model parameters (e.g., the sample median),
the JK method may produce an inconsistent estimator of the 
MSPE (more precisely, of the variance type  term $\mtoi$), while  the bootstrap
 based LM1 and LM2 methods 
would still work (cf. Ghosh et al. (1984)). From this point of view,
 the proposed  method of MSPE estimation has a  wider range of 
validity than the JK method.

%A second important observation is that  in contrast to  
%the PR- and the JK  methods which may produce a negative estimate of the 
%MSPE (although not in the normal case studied by Prasad and Rao (1990) 
%or Datta and Lahiri (2000)), the estimates of the MSPE produced by the 
%proposed method is always nonnegative. This is  important for  
% constructing prediction  intervals, where a zero or a negative 
%estimate of the MSPE is useless. However,  the 
%approximation based LM2 method  is computationally 
%more extensive than  the JK method (and the PR method, when it is
%applicable). 

In this paper, we also %establish 
%investigate theoretical properties of the LM2
%method and 
prove  that the  proposed  estimator of the MSPE attains
 the same level of  asymptotic  accuracy as the existing methods in correcting 
the bias of  the simple  MSPE estimator.  We also report the results of 
a small simulation study and provide some guidelines
for implementing the methodology in practice. 
%(The codes 
%for applying the LM2 method can be obtained from one of the 
%authors by email {\it taps\@iastate.edu}.) 
 In summary, the proposed 
method allows a user to routinely derive second order accurate, nonnegative 
estimates of the MSPE in small area estimation problems, 
without requiring  any analytical work on the part of the user.\\[.1in]
\centerline {\sc Acknowledgement} 
%This work is  patially supported by a 
%contract from the National Center for Health Statistics.
The authors thank three referees for their 
constructive criticism that led to a vast improvement of an
earlier draft of the paper. The authors also thank 
 Douglas  Williams  for some  helpful 
discussions  and for facilitating 
 the project.

\vspace*{.1in}

\begin{center}
{\large{\sc References}}
\end{center}

\begin{description}
\item 
{\sc Booth}, J.G., \& {\sc Hobert}, J.P. (1998). Standard errors of prediction in 
generalized linear mixed models. {\it J. Am. Statist. Assoc.} {\bf 83}, 28-36.

%\item
%{\sc Breslow}, N. E., \&  {\sc Clayton}, D. G.   (1993). Approximate inference 
%in generalized 
%linear mixed models. {\it J. Am. Statist. Assoc. }{\bf 88}, 9-25.

\item
{\sc Butar}, F.B., \& {\sc Lahiri}, P. (2003). On measures of uncertainty of 
empirical Bayes small area estimators. {\it J. Statist. 
Plan.  Infer.} {\bf 112},   63-76.

\item
{\sc Datta}, G.S. \&  {\sc Lahiri, P.} (2000). A unified measure of
uncertainty of estimated best linear unbiased predictors in small area
estimation problems. {\it  Statist. Sinica}, {\bf 10}, 623-27.

\item 
{\sc Datta}, G.S., {\sc Rao}, J.N.K., \& {\sc Smith}, D.D. (2005).
 On measuring the 
variability of small area estimators under a basic area level model. {\it Biometrika}
 {\bf 92}, 183-96.
 
\item
{\sc Efron}, B. (1978). Regression and ANOVA with zero-one data: Measures 
of residual variation. {\it J. Am. Statist. Assoc.}  {\bf 73}, 
113-21.

\item
{\sc Efron}, B.   (1979). Bootstrap methods: Another look at the jackknife.
{\it  Ann. Statist.}  {\bf 7},   1-26.

\item
{\sc Efron}, B. (1986). Double exponential families and their use in 
generalized linear regression. {\it J. Am. Statist. Assoc.}  {\bf 81}, 709-21.

\item
{\sc Efron}, B.  \&  {\sc Tibshirani}, R.    (1993).
{\it An introduction to the bootstrap.} Chapman \& Hall Ltd., New York. 

\item
{\sc Fay}, R. E. \&  {\sc Herriot}, R. A. (1979). Estimates of income for small 
places: An application of James-Stein procedures to census data.
 {\it J. Am. Statist. Assoc.} {\bf  74},   269-77.

\item
{\sc Ghosh, M.}, \& {\sc Maiti, T.} (2004). Small-area estimation based on natural 
exponential family quadratic variance function models and survey weights. 
{\it Biometrika} {\bf 91}, 95-112.

\item
{\sc Ghosh}, M.,  {\sc Parr}, W. C.,  {\sc Singh}, K., \&  {\sc Babu}, G.J.  (1984).
A note on bootstrapping the sample median. {\it  Ann. Statist.}
{\bf   12},   1130-35.
 
%\item
%{\sc Ghosh}, M.,  {\sc Natarajan}, K.,  {\sc Stroud}, T. W. F., \&  {\sc Carlin}, %B. P.   (1998).
%Generalized linear models for small-area estimation. 
%{\it J. Am. Statist. Assoc.}  {\bf 93}   273-82.

\item
{\sc H\"{a}rdle}, W. (1991). {\it Smoothing Techniques: With implementation 
in S}. Springer, New York, NY.

\item
{\sc Jiang}, J., {\sc Lahiri}, P., \&  {\sc Wan}, S-M. (2002). A unified 
jackknife theory for empirical best prediction with M-estimation. 
{\it  Ann.   Statist.} {\bf 30}, 1782-810.

 \item
{\sc Lahiri}, S.N. \& {\sc Maiti}, T. (2003). Nonnegative Mean Squared Error 
prediction. {\it Preprint}. Posted at
{\it  http://arxiv.org/abs/math.ST/0604075}.

\item
{\sc Lahiri}, S.N., {\sc Maiti, T.}, {\sc Katzoff, M.}, \& {\sc Parsons,V.} 
(2006). Resampling based empirical prediction: An application to small area 
estimation. Posted at \\
{\it  http://arxiv.org/abs/math.ST/0604513}.

%\item
%  {\sc Lee}, Y.  , \&  {\sc Nelder}, J. A.   (1996). Hierarchical generalized 
%linear models (with discussion).
%{\it  J. R. Statist. Soc.  B}  {\bf 58},   619-56.

\item
{\sc McCulloch}, C.E., \& {\sc Searle}, S.R. (2001). 
{\em Generalized, Linear, and Mixed Models}
 New York: Wiley.

\item 
{\sc Pfeffermann}, D. \& {\sc Glickman}, H. (2004). Mean squared error approximation in 
small area estimation by use of parametric and nonparametric bootstrap. 
{\it Proc. Sec.  Survey Res. Meth., Am. 
Statist. Assoc.}

\item 
{\sc Pefferemann}, D., \& {\sc Tiller}, R.B. (2005). Bootstrap approximation to prediction 
MSE for state-space models with estimated parameters. {\it J. Time Ser. Analysis},
 {\bf 26}, 893-16.

\item
{\sc Prasad}, N.G.N. \&  {\sc Rao}, J.N.K. (1990). The estimation of the 
mean squared error of small area estimators. {\it  J. Am. Statist. Assoc.} 
{\bf 68},  67-72. 

\item
{\sc Rao, J.N.K.} (2003). {\it Small Area Estimation}.  Wiley, New York.

\item {\sc Slud, E.V.}  (2000). Comparison of aggregate versus unit-level 
models for small area estimation. {\it Proc. Sec.  Survey Res. Meth., Am. 
Statist. Assoc.}

\item 
{\sc Slud}, E.V., \&  {\sc Maiti}, T. (2006). MSE estimation in transformed 
Fay-Herriot models. {\it J. R. 
Statist. Soc.  B} {\bf 68}, 239-57.

\item
{\sc You, Y.} \& {\sc Rao, J.N.K.} (2002). Small area estimation using unmatched 
sampling and linking models. {\it Can. J. Statist.} {\bf 30}, 3-15.
\end{description}

% \end{document}

%%%%%%%%%%%%%%%%%%%%%%%%%%%%%%%%%%%%%%%%%%%%%%%%%%%%%%%%%%%%%%
\vspace*{.1in}

\begin{center}
{\Large{\sc Appendix A: Proofs}}\\[.2in]
\end{center}
\renewcommand{\theequation}{A.\arabic{equation}}
\setcounter{equation}{0}

Let $\IN=\{1,2,\ldots\}$.
% denote the set of all positive integers.
In the proofs,  we suppress  dependence of various quantities
on $m$ unless there is a chance of confusion and write  
$C, C(\cdot)$ to  denote generic  positive constants
 that depend on their arguments (if any), but not on
$i\in\{1,\ldots,m\}$ or $m$.
\\

\noindent
{\bf Lemma 1:} Let $X_1,\ldots,X_n$ ($n\geq 1$)  be a collection of 
iid random variables with $E|X_1|^{1+\eta}<\infty$  for some 
$\eta\in (0,1]$. Let  $\mu = EX_1$, $\bar{X}_n = n^{-1}\sum_{i=1}^n X_i$
and $\rho = (E|X_1|^{1+\eta})^{\frac{1}{1+\eta}}$. Then 
\be
%%%%%%%%%%%%%%%%%%%%%%%%%%%%%%%%%%%%%%%%%%%%%%%%%%%%%%%%%%%%%%%%%%% eqn (A.1)
E|\bar{X}_n -\mu| \leq 3 \rho\, n ^{{\eta}/{(1+\eta)}}
\qmq{for all} n\geq 1.
\ee

\noindent
{\bf Proof:} If $\rho =0$, then 
%$P(X_1=0)=1$ and hence 
%$
%P(\bar{X}_n=0)=1$. In this case, 
(A.1) holds trivially. 
Hence, suppose that $\rho> 0$. With 
$c_n = \rho n^{{1}/{(1+\eta)}}$, let $X_{1i}= X_iI(|X_i|\leq c_n)$,
 $X_{2i}= X_i -X_{1i}$, $1\leq i\leq n$, and 
 $\bar{W}_{kn}
=n^{-1}\sum_{i=1}^n (X_{ki}- E X_{ki})$, $k=1,2$. Then, 
%$\bar{X}_n =\bar{W}_{1n}+\bar{W}_{2n}$ and 
$
E|\bar{X}_n -\mu| \leq 
	E|\bar{W}_{1n}|+E|\bar{W}_{2n}|
%\\
%&\leq& 
%	(E|\bar{W}_{1n}|^2)^{1/2} + E|\bar{W}_{2n}|\\
\leq 	
	(n^{-1}E|X_{11}|^2)^{1/2} + 2E|X_{21}|
\leq 
	(n^{-1}c_n^{1-\eta}E|X_{1}|^{1+\eta})^{1/2}
 + 2E|X_{1}|^{1+\eta} c_n^{-\eta}
\leq 3 \rho\, n ^{{\eta}/{(1+\eta)}}.
$
\\

\noindent
{\bf Lemma 2:} For random vectors $\uX$ and $\uY$
%be a $r$- and 
% $s$-dimensional ($r,s\in\IN$) random vectors
on a common probability space with 
 $E|g(\uY)|^{\alpha} <\infty$ for some $g:\Re^s\raw \Re$ 
and   $\alpha\in [1,\infty)$,~
$
E| E\{g(\uY)|\uX\} - g(\uY)|^{\alpha}
\leq 2^{\alpha} E|g(\uY)|^{\alpha}.
$
\\[.1in]
\noindent
{\bf Proof:}
 Follows from H\"older's and 
 conditional 
Jensen's inequalities.
\\

\noindent
{\bf Proof of Theorem 1:} ~(i) By (C.1)(ii) and Taylor's expansion, 
%the smoothness conditions  on 
%$M_{1i}(\cdot)$, 
for some $u_{1i}, u_{2i}\in [-1,1]$,  
\begin{eqnarray}
%%%%%%%%%%%%%%%%%%%%%%%%%%%%%%%%%%%%%%%%%%%%%%%%%%%%%%%%%%%%%%%%%%%%%%%%%%%%%% (A.2) - r2
&&E\Big | \Big\{ M_{1i}(\hud+ \z \ue_j) -  M_{1i}(\hud- \z \ue_j) 
		\Big\} -2\z M_{1i}^{(j)}(\hud)\Big |\nn\\
%&=&	
%	E\Big| \Big[ \Big\{ M_{1i}(\hud)+\z M_{1i}^{(j)}(\hud)
%		+\frac{[\z]^2}{2} M_{1i}^{(j,j)}(\hud+u_{1i}\z\ue_j)\Big\}\nn\\
%&&
%	- \Big\{ M_{1i}(\hud)-\z M_{1i}^{(j)}(\hud)
%		+\frac{[\z]^2}{2} M_{1i}^{(j,j)}(\hud+u_{2i}\z\ue_j)\Big\}\Big]
%			-2\z M_{1i}^{(j)}(\hud)\Big|\nn\\
&=&	
	2^{-1}{[\z]^2} E\Big|  M_{1i}^{(j,j)}(\hud+u_{1i}\z\ue_j)
		-   M_{1i}^{(j,j)}(\hud+u_{2i}\z\ue_j)\Big|
\leq
	C_1 [\z]^2 
\end{eqnarray}
for all $i=1,\ldots,m$, $ m\geq 1$. Since
$
E_{\udelta_0} M_{1i}^*(\udelta_0) = M_{1i}(\udelta_0)
$
for all $\udelta_0\in\Delta$, by Lemmas 1 
and 2,
\begin{eqnarray}
%%%%%%%%%%%%%%%%%%%%%%%%%%%%%%%%%%%%%%%%%%%%%%%%%%%%%%%%%%%%%%%%%%%%%%%%%%  (A.3) - r2
&&
	E\Big| M_{1i}^*(\hud_1) -M_{1i}(\hud_1)\Big|
%\nn\\
%&=&
%	E\left[ E_{\hud_1} \Big | 
%	 	\avelno \{\xi_i(y_i^{*l}) - h(\theta_i^{*l})\}^2 
%			- M_{1i}(\hud_1)\Big|\right]\nn\\
\leq  
	3 E\{ E_{\hud_1} | \xi_i(y_i^{*1}) 
				- h(\theta_i^{*1})|^{2+2\eta} 
					N_0^{-{\eta}/{(1+\eta)}}\}\nn\\
&\leq & 
	C(\eta)	 E\{   E_{\hud_1}|h(\theta_i^{*1})|^{2+2\eta} 
					N_0^{-{\eta}/{(1+\eta)}}\}
\leq 
	C(\eta) N_0^{-{\eta}/{(1+\eta)}},
\end{eqnarray}
where $\hud_{1} \equiv \hud +\z \ue_{j}$.
Using (A.3) and  similar arguments for $M_{1i}^*(\hud -\z \ue_{j})$, 
we get
\be
%%%%%%%%%%%%%%%%%%%%%%%%%%%%%%%%%%%%%%%%%%%%%%%%%%%%%%%%%%%%%%%%%%%%%%%%  (A.4) - r2
E|M_{1i}^{(j)*} (\hud)
	 -(2\z)^{-1}\{ M_{1i}(\hud+\z\ue_j)
			 -  M_{1i}(\hud-\z\ue_j)\}|
	\leq (2\z)^{-1}\{C(\eta) N_0^{-{\eta}/{(1+\eta)}}\}
\ee
uniformly in $i=1,\ldots,m$, $m\geq 1$.
Part (i) of the theorem now follows from (A.2)-(A.4).

Next consider (ii).  By arguments similar to (A.2),~
$E | \{ M_{1i}(\hud+ \z \ue_j) + M_{1i}(\hud- \z \ue_j) 
	 -2 M_{1i}(\hud)	\}
		-(\z)^2 M_{1i}^{(j,j)}(\hud)|
\leq
	C\cdot  (\z)^3
$
uniformly in $i=1,\ldots,m$, $m\geq 1$. Also, using Lemma 1, the linearity
of $M_{1i}^{(j,j)*}(\hud)$ in $M_{1i}^{(j)*}(\cdot)$ and arguments similar to 
(A.3), one can
show that 
\be
%%%%%%%%%%%%%%%%%%%%%%%%%%%%%%%%%%%%%%%%%%%%%%%%%%%%%%%%%%%%%%%%%%%%%%%%  (A.5) - r2
E |M_{1i}^{(j,j)*}(\hud) -(\z)^{-2}
\{ M_{1i}(\hud+ \z \ue_j) + M_{1i}(\hud- \z \ue_j) 
	 -2\z M_{1i}^{(j)}(\hud)	\}|
	\leq C(\eta) (\z)^{-2}  N_0^{-{\eta}/{(1+\eta)}}
\ee
uniformly in $i=1,\ldots,m$, $m\geq 1$. Hence, part (ii) holds for all
$j,r\in\{1,\ldots,k\}$ with $j=r$. Next fix $1\leq j\neq r \leq k$. 
Define
$
M_{1i}^{(j,r)\dagger}(\hat{\udelta}) 
= (2[\z]^{2})^{-1}[
	\{
		 M_{1i}(\hat{\udelta} + \z\ue_{j,r})
	+ M_{1i}(\hat{\udelta} - \z\ue_{j,r})
	- 2M_{1i}(\hat{\udelta})
		\}
		- [\z]^{2}\{
		M_{1i}^{(j,j)}(\hat{\udelta}) 
		+ M_{1i}^{(r,r)}(\hat{\udelta}) 
				\}
	]$. 
By 
Taylor's expansion 
% it is easy to show that 
\be
%%%%%%%%%%%%%%%%%%%%%%%%%%%%%%%%%%%%%%%%%%%%%%%%%%%%%%%%%%%%%%%%%%%%%%%%  (A.6) - r2
E| M_{1i}^{(j,r)\dagger}(\hat{\udelta}) 
	 - M_{1i}^{(j,r)}(\hat{\udelta}) | (2[\z]^2)
\leq  C[\z]^3.
\ee
Now  using (A.6) and arguments similar to (A.5), one can complete the 
proof of (ii).
% for the $j\neq r$ case. We omit the 
%routine details.
\\

\noindent
{\bf Proof of Theorem 2:}
 Note that $\ecd (\bd^{*1}) -\hud
=\ub(\hud) = \hub$ and  $\ecd (V^{*1}) =\hV$.
 Hence,  for any $j\in\{1,\ldots,k\}$,
%\begin{eqnarray*}
$
E| b^*(j) -\hat{b}(j)|^2 
\leq 
	N_0^{-1}E\{\ecd( \delta^{*l}(j) 
		-\hat{\delta}(j))^2\}
\leq
	2 N_0^{-1}E\{\rho_m(\hud;2) +\|\hud\|^2\}
=
	 O(N_0^{-1})
$.
Similarly,
 by Lemma 1, $E|V^*(j,r) -\hV(j,r)|$ is bounded above by
\begin{eqnarray*}
\vspace*{-.2in}
&&
E\Big[ \ecd |\avelno \delta^{*l}(j)\delta^{*l}(r)
		- \hat{\delta}(j) \hat{\delta}(r)|\\
&&\hspace{.1in}
+ \{\ecd |\bar{d}_j^*
			-\hat{\delta}(j)|^2\}^{1/2}
		\{\ecd (\bar{d}_r^*)^2\}^{1/2}
	+ |\hat{\delta}(j)|( \ecd | \bar{d}_r^*
			-\hat{\delta}(r)|^2)^{1/2}\Big]\\
&\leq&
	C(\eta)\Big[ E\Big\{\rho_m(\hud;2+2\eta) 
				+ \rho_m(\bd;2+2\eta)\Big\} 
					N_0^{-\frac{\eta}{1+\eta}}
	+ E \Big\{ \rho_m(\hud;2) 
				+ \rho_m(\bd;2)\Big\} 
					N_0^{-1/2}\Big],
\end{eqnarray*}
for any $j,r\in\{1,\ldots,k\}$, where 
$\bar{d}_j^* = \avelno\delta^{*l}(j)$.  
Theorem 2   follows from these bounds.\\

\noindent
{\bf Lemma 3 :}
Suppose that
condition (C.3) holds. Then, for any  $\gamma 
\in (0,\eta)$,
$
E\|\hub -\ub\|^{1+\gamma} = o(m^{-(1+\gamma)})$ 
{and}
 $E\|\hV -V\|^{1+\gamma} = o(m^{-(1+\gamma)})$.
\\[.1in]
\noindent
{\bf Proof:} Fix $\gamma\in(0,\eta)$. Note that 
$mE\|\hud -\bd\|^2 \leq C[\|\ua_m(\bd)\|^2 + \|\Sigma_m(\bd)\|^2]
\leq C G_2(\bd)<\infty$. Hence, $\hud\raw \bd$ in mean sqrare and 
therefore,  by the equicontinuity condition,
$\|\ua_m(\hud)- \ua_m(\bd)\|$ and $\|\Sigma_m(\hud)
-\Sigma_m(\bd)\|$ both converge to zero in probability
under $\bd$. Further, the sequence $\{G_2(\hud)^{1+\gamma}\}$
is uniformly integrable. Hence, by the (extended) Dominated Convergence
Theorem,~
$
[E\|\ua_m(\hud)-\ua_m(\bd)\|^{1+\gamma} + 
E\|\Sigma_m(\hud) -\Sigma_m(\bd)\|^{1+\gamma}
] ~\raw 0$ as  $\mti
$, proving the lemma.\\

\noindent
{\bf Proof of Theorem 3:} First we show that 
\be
%%%%%%%%%%%%%%%%%%%%%%%%%%%%%%%%%%%%%%%%%%%%%%%%%%%%%%%%%%%%%%%%%%%%%%%%%%% eqn (A.7) - r2
\max_{1\leq i\leq m} E| M_{1i}^*(\cdis) - M_{1i}(\cdi)| 
+
\max_{1\leq i\leq m} E | M_{2i}^*(\hud) - M_{2i}(\hud)| 
=
	o(m^{-1}).
\ee
Consider the  first term on the left side.  By arguments similar to (A.5), 
$
\max_{1\leq i\leq m} E| M_{1i}^*(\cdis) - M_{1i}(\cdis)|
\leq C(\eta) N_0^{-{\eta}/{(1+\eta)}}
$.
Next, write $\ais=\{ \bdis\in\De\}\cap \{|\moi^{(s)*}(\hud)|^{-1}
\leq (1+\log m)^2\}$ and $\ai =\{ \bdi\in\De\}\cap \{|\moi^{(s)}(\hud)|^{-1}
\leq (1+\log m)^2\}$, $1\leq i\leq m$, $m\geq 1$. Then, using (4.13), it can be shown that 
\bea
%%%%%%%%%%%%%%%%%%%%%%%%%%%%%%%%%%%%%%%%%%%%%%%%%%%%%%%%%%%%%%%%%%%%%%%%% eqn (A.8) - r2
E| \moi(\cdis) - \moi(\cdi)|
&\leq&
	E| \moi(\bdis) -\moi(\bdi)| I(\ais\cap\ai)
		+ E[ \moi(\hud)\{I(\ai^c) + I([\ais]^c)\}]\nn\\
&&
	+ E\moi(\bdi)I([\ais]^c\cap\ai) + E\moi(\bdis)I(\ais\cap\ai^c)\nn\\
&\equiv&
	R_{1i}+	R_{2i}+	R_{3i}+	R_{4i},\qmq{say.}
\eea
By (C.1), (C.2)
(with $\eta=1$), (C.3) and  arguments similar to the proof of Theorem 1, 
one gets ~
${\max}_{1\leq i\leq m} E| \moi^{(j)*}(\hud) - \moi^{(j)}(\hud)|^2
= O([\z]^2 + [\z]^{-2}N_0^{-1})$, 
${\max}_{1\leq i\leq m} E| \moi^{(j,j)*}(\hud) - \moi^{(j,j)}(\hud)|^2
= O([\z]^2 + [\z]^{-4}N_0^{-1})$,
and 
$
E\|\ub\|^2 + E\|V\|^2 = O(m^{-2})
$. Now using  the above bounds, it can be shown (cf. (A.17), Lahiri et al. (2006)) that 
\bea
\mbox{max}_{1\leq i\leq m} R_{1i} 
\leq 
%%%%%%%%%%%%%%%%%%%%%%%%%%%%%%%%%%%%%%%%%%%%%%%%%%%%%%%%%%%%%%%%%%%%%%%%% eqn (A.9) - r2
	C_1 \mbox{max}_{1\leq i\leq m} E\|\bdis -\bdi\|I(\ais\cap \ai)= o(m^{-1}).
\eea
Since $|\moi^{(s)}(\bd)| >C_0$, there exist $\ep_1,\ep_2\in(0,\infty)$
such that $|\moi^{(s)}(\ux)|>\ep_1$ for all $\ux\in\De$ with 
$\|\ux-\bd\|\leq \ep_2$. Hence, by (C.1),
there exists a $C = C(\ep_1)\in(0,\infty)$ such that
 on the set 
$\{\|\hud -\bd\| \leq \ep_2\}$,\,
$
\|\bdi -\bd\|\leq C [\|\hub\|+\|\hV\|]
$
for all $i=1,\ldots,m$, $m\geq 1$. Hence, for any $\ep >0$, by (C.1) and (C.4),
(cf. (A.18)-(A.19), Lahiri, et al. (2006))
\bea
%%%%%%%%%%%%%%%%%%%%%%%%%%%%%%%%%%%%%%%%%%%%%%%%%%%%%%%%%%%%%%%%%%%%%%%% eqn (A.10) - r2
&&\max_{1\leq i\leq m} P(\|\bdi -\bd\|>\ep)
\leq
	P(\|\hud -\bd\|>\ep_2) + P(C[\|\hub\|+\|\hV\|] >\ep)
= O(m^{-(1+\eta)}),\quad\quad\\
%%%%%%%%%%%%%%%%%%%%%%%%%%%%%%%%%%%%%%%%%%%%%%%%%%%%%%%%%%%%%%%%%%%%%%%% eqn (A.11) - r2
&&\max_{1\leq i\leq m} P\Big(|\moi^{(s)}(\hud)|\leq (1+\log m)^{-2}\Big)
%\leq 
%	\max_{1\leq i\leq m} P\Big(|\moi^{(s)}(\hud) - \moi^{(s)}(\bd)| 
%			>\ep_1/2\Big)
\leq 
	2P\Big(\|\hud -\bd\|> C(\ep_1,\ep_2)\Big)
		= O\Big(m^{-(1+\eta)}\Big).\quad\quad
\eea
Hence, it follows that 
\be
%%%%%%%%%%%%%%%%%%%%%%%%%%%%%%%%%%%%%%%%%%%%%%%%%%%%%%%%%%%%%%%%%%%%%%%% eqn (A.12) - r2
\max_{1\leq i\leq m} P(\ai^c) =  O\Big(m^{-(1+\eta)}\Big).
\ee
We now obtain a similar bound on $P([\ais]^c)$. Since $\bd$ is an 
interior point of $\De$, there exists a $\ep_3\in(0,\infty)$ such that 
$\{\ux: \|\ux -\bd\|  \leq \ep_3\}\subset\De$. Let $A_{1i}^*
=\{\bdis\in\De\}$ and  $A_{2i}^*
=\{|\moi^{(s)*}(\hud)|^{-1} \leq (1+\log m)^2\}$. By (5.3), and (A.7)-(A.12),
uniformly over  $i=1,\ldots,m$,
\bea
%%%%%%%%%%%%%%%%%%%%%%%%%%%%%%%%%%%%%%%%%%%%%%%%%%%%%%%%%%%%%%%%% eqn (A.13) - r2
&&P([\ais]^c) 
\leq	
	P(A_{1i}^{*c}\cap A_{2i}^{*}\cap \ai) + 
		+P(\ai^c) + P(A_{2i}^{*c})\nn\\
&\leq&
	P(\|\bdi -\bd\| >{\ep_3}/{2}) 
		+ {2}{\ep_3}^{-1}E\|\bdis -\bdi\|I(A_{2i}^*\cap\ai)
				+P(\ai^c)\nn\\
&&
			+\Big[ P\Big(|\moi^{(s)*}(\hud)
				- \moi^{(s)}(\hud)|
			> \frac{C_0}{2}-\frac{1}{(1+\log m)^{2}}	\Big)
			+P\Big(|\moi^{(s)}(\hud) -\moi^{(s)}(\bd)|
				>\frac{C_0}{2}\Big)
				\Big]\nn\\
& =& O\Big(m^{-(1+\eta)} 
			+ (\log m)^4[\z+(\z N_0^{1/2})^{-1}]\Big).
\eea
Now using   (A.12), (A.13) and  condition (C.1), with $a_i^2\equiv P(\ai^c)+P([\ais]^c)$,
we have 
\bea
%%%%%%%%%%%%%%%%%%%%%%%%%%%%%%%%%%%%%%%%%%%%%%%%%%%%%%%%%%%%%%%%% eqn (A.14) - r2
R_{2i}
&\leq&
	C \{a_i^2 % P(\ai^c)+P([\ais]^c)
	+ C_1 a_i%\{P(\ai^c)+P([\ais]^c)\}^{frac12}
(E\|\hud-\bd\|^2)^{1/2}
\}
= o(m^{-1}),\\
R_{3i}
&\leq&
	E\Big|\moi(\bdi) -
		\moi(\hud)\Big|I([\ais]^c\cap \ai) + 
			R_{2i}
= o(m^{-1}),\\
 R_{4i}
&\leq&
E|\moi(\bdis) -
		\moi(\hud)|I(\ais\cap \ai^c) + 
			R_{2i} = o(m^{-1}),
\eea
uniformly in $i\in\{1,\ldots,n\}$ (cf. (A.23)-(A.25), Lahiri et al. (2006)). 
By (A.8), (A.9), and (A.14)-(A.16), the first term on the left of (A.7) is $o(m^{-1})$.
The upper  bound on the  other term  
on the left of (A.7)
follows from condition  (C.1),   
the independence of the resampled vectors $(y_1^{*l},\ldots,
y_m^{*l})$ for $l=1,\ldots, N_0$ and the fact 
$
\ecd( \xi_i(y_i^{*1};\udelta^{*1}) - \xi_i(y_i^{*1};\hud))^2
= \mtoi(\hud)
$.
Hence (A.7) is proved  which, in turn, implies that 
$\max_{1\leq i\leq m} E \Big|mspe_{\mbox{{\sc new}}}(\hat{\beta}_i)
         - mspe_{\mbox{{\sc lm}}}(\hat{\beta}_i)\Big|
=o(m^{-1})
$.
Next  define 
the preliminary titled estimator 
$\bdi$  for the LM method by using 
 the bias and the 
variance estimators  $\hub =\ub(\hud)$ and 
$\hV= V(\hud)$. Note that with this choice
of $\hub$ and $\hV$, the regularity conditions 
for the validity of  Theorem 3 of LM follow from
conditions (C.1)-(C.4)  and Lemma 3 above. Hence, (5.4) 
follows  from Theorem 3
of LM.

\newpage
{\bf Appendix B}

In this section, the simulation results are presented into subclasses as per the 
request of a referee. For example, in model I and Model III, the small areas are 
grouped into 3 classes having eaual sampling variances, denoted as G1, G2 and G3.
Thus each group represent 5 areas and summary results are presented for each group. 
But for model II, 3 representative areas are chosen, namely the areas for $n_i=6$, 
$n_i=16$ and $n_i=36$. Though they are not group in a true sense, they are also represented as G1, G2 and G3 in the tables for convenience. Note that, in this case the estimates represent only thsese selected three areas, not the averages.

The Table 1b represnts the simulated bias and MSPE. For model I, the third group 
has higher bias and vice versa for model III. For model II, G3, the highest sample size has lowest bias. Interms of MSPE, for all the models, G1 is the highest, althogh the results between the goups are not drastically different. Also the kernel based method and the closed form formulas (wherever applicable) performs equally.

\begin{center}
 Table 1b.
{\it  Absolute relative bias ($T_1$) and empirical MSPE 
($T_2$) for the EBP. Results using the 
   kernel based approximations (3.6) and (3.7)   are reported within the parentheses.}
% under the
% results using standard methods.

\begin{tabular}{ccccccc}\\ 
       & &\multicolumn{2}{c}{Model I}& Model II 
& \multicolumn{2}{c}{Model III} \\ 
Measures & Group & Median & Mean   &  & Median & Mean\\ 

     &G1 & 2.201  & 2.087 & 0.276  & --- & --- \\
     &   &(2.119) &(2.005)&(0.271) &(1.821)&(1.194) \\ 
$T_1$&G2 & 1.804  & 2.196 & 0.197  & --- & --- \\
     &   &(2.030) &(2.066)&(0.199) &(1.001)&(1.034) \\
     &G3 & 2.476  & 4.846 & 0.156  & ---  &  ---  \\  
     &   &(2.631) &(4.265)&(0.155) &(.840)& (0.825) \\
     &   &        &       &        &      &         \\
     &G1 & 0.456  & 0.435 & 0.015  & ---  & ---\\
     &   &(0.468) &(0.483)&(0.019) &(0.300)& (0.298)\\
$T_2$&G2 & 0.372  & 0.360 & 0.013  & ---   & --- \\
     &   &(0.375) &(0.362)&(0.012) &(0.272) & (0.282) \\
     &G3 & 0.234  & 0.240 & 0.003  & --- & --- \\
     &   &(0.244) &(0.243)&(0.003) &(0.250) & (0.245) \\ 
\end{tabular}
\end{center}

\newpage

The relative bias and the coefficient of variations of the MSPE estimates are 
presented in Table 2b. The results for LM1 and LM2 are based on Remark 3 modification. However, they are fairly close when (4.6) was used instead. For all the 
groups the JLW shows slightly higher bias and CV compared to others. LM1 and PR/DL 
performs equally well both in terms of bias and CV, LM2 has little higher CV for model I. For model II, CV under JLW is higher than that under LM2. For model III, LM2 performs well for all the groups. For large sample size, the CV under JLW is small yet larger than other methods.

\begin{center}
 Table 2b.
{\it  Relative biases ($T_3$) and coefficient of variations ($T_4$)
 for the bias corrected estimators of the MSPE.
Entries for LM1 and  LM2 are based on Remark 3 modification.
}

\begin{tabular}{cccccccc}\\ 
       &        &   &\multicolumn{2}{c}{Model I}& Model II 
& \multicolumn{2}{c}{Model III} \\ 
Method & Measures & Group & Median & Mean &    & Median & Mean \\ 

PR/DL  & $T_3$ &G1& 0.016  & 0.090 &  ---  & --- &---\\
       &       &G2& 0.008  & 0.063 &  ---  & --- & --- \\
       &       &G3& 0.106  & 0.084 & ---   & --- & --- \\
       & $T_4$ &G1& 0.184  & 0.252 &  ---  & --- & ---\\
       &       &G2& 0.151  & 0.203 & ---   & --- & --- \\
       &       &G3& 0.119  & 0.113 & ---   & --- & --- \\

JK     & $T_3$ &G1& 0.287  & 0.243 & -0.190 & --- & --- \\
       &       &G2& 0.124  & 0.152 & -0.083 & --- & --- \\
       &       &G3& 0.124  & 0.173 &  0.025 & --- & --- \\

       & $T_4$ &G1& 0.924  & 0.705 &  1.532 & --- & --- \\
       &       &G2& 0.379  & 0.449 &  0.752 & --- & --- \\
       &       &G3& 0.366  & 0.419 &  0.360 & --- & --- \\

LM1    & $T_3$ &G1& -0.000 & 0.072 &  ---  & --- & ---\\
       &       &G2& -0.017 & 0.036 & ---   & --- & --- \\
       &       &G3&  0.061 & 0.041 & ---   & --- & --- \\

       & $T_4$ &G1&  0.190 & 0.246 &  --- & --- & ---\\
       &       &G2&  0.163 & 0.196 &  --- & --- & --- \\
       &       &G3&  0.083 & 0.095 & ---  & --- & --- \\

LM2    & $T_3$ &G1& -0.093 &-0.018 & -0.148 & -0.005& -0.000  \\
       &       &G2& -0.102 &-0.039 &  0.094 &  0.102&  0.009 \\
       &       &G3&  0.003 &-0.016 &  0.054 &  0.152&  0.108 \\
 
       & $T_4$ &G1&  0.276 &  0.300 &  0.154&  0.414&  0.368 \\
% \hspace{0in}0.319 (0.310)& \hspace{0in} 0.368 (0.298) \\
%       &          &        &          &       &      & 0.310& 0.298\\
        &      &G2&  0.263 &  0.274 &  0.746 & 0.102 & 0.009 \\
% \hspace{0in}0.319 (0.310)& \hspace{0in} 0.368 (0.298) \\
        &      &G3&  0.201 &  0.202 &  0.054 &  0.309 & 0.202 \\ 
% \hspace{0in}0.319 (0.310)& \hspace{0in} 0.368 (0.298) \\

\end{tabular}
\end{center}

\end{document}